\documentclass[12pt]{amsart}
\usepackage{amsmath,amsfonts,amscd,amssymb,latexsym,amsthm}

\title{Algebraic Geometry over Lie Algebras}
\author{Ilya V. Kazachkov}
\address{I. V. Kazachkov\newline
Department of Mathematics and Statistics \newline McGill University
\newline
805 Sherbrooke st. West\newline Montreal, Quebec, H3A 2K6, Canada}
\email{ilya.kazachkov@gmail.com}
\thanks{Supported by a grant from the London Mathematical Society.}

\newtheorem{lem}{Lemma}[section]
\newtheorem{thm}[lem]{Theorem}
\newtheorem{cor}[lem]{Corollary}
\newtheorem{prop}[lem]{Proposition}

\theoremstyle{definition}

\newtheorem{defn}{Definition}[section]
\newtheorem{expl}{Example}[section]
\newtheorem{rem}{Remark}[section]

\newtheorem*{problem}{Problem}

\newcommand{\x}{\xi}

\newcommand{\D}{\mathcal{D}}
\newcommand{\F}{\mathfrak{F}}
\newcommand{\V}{\ensuremath{\mathbb{V}}}
\newcommand{\FrGr}{\ensuremath{\mathbb{F}}}

\renewcommand{\hom}{\texttt{Hom}}
\newcommand{\rad}{\texttt{Rad}}
\newcommand{\AS}{\texttt{AS}}

\newcommand{\CA}{\texttt{CA}}
\newcommand{\Ff}{\texttt{F}}
\newcommand{\Gg}{\texttt{G}}
\newcommand{\RCA}{\texttt{CPA}}

\newcommand{\K}{\mathcal{K}}
\newcommand{\im}{\hbox{im}}

\newcommand{\AX}{A \left[ X \right]}

\newcommand{\Fx}{F \left[ x \right]}
\newcommand{\pvar}{\texttt{-pvar}}
\newcommand{\qvar}{\texttt{-qvar}}
\newcommand{\var}{\texttt{-var}}
\newcommand{\sFit}{\textsf{Fit}}
\newcommand{\Fit}{\texttt{Fit}}
\newcommand{\id}{\texttt{id}}
\newcommand{\ucl}{\texttt{ucl}}
\newcommand{\fucl}{\texttt{-ucl}}

\newcommand{\Q}[1]{\marginpar{\textbf{$\Leftrightarrow$}}\ }
\newcommand{\bi}{\begin{itemize}}
\newcommand{\ei}{\end{itemize}}
\begin{document}

\maketitle

\section*{Introduction}
\label{chap:I}

What is algebraic geometry over algebraic systems? Many important
relations between elements of a given algebraic system
$\mathcal{A}$ can be expressed by systems of equations over
$\mathcal{A}$. The solution sets of such systems are called
\emph{algebraic sets} over $\mathcal{A}$. Algebraic sets  over
$\mathcal{A}$ form a category, if we take for morphisms polynomial
functions in the sense of Definition \ref{defn:I-polmap} below. As
a discipline, algebraic geometry over $\mathcal{A}$ studies
structural properties of this category.  The principle example is,
of course, algebraic geometry over fields. The foundations of
algebraic geometry over groups were laid by Baumslag, Myasnikov
and Remeslennikov \cite{alg1,alg2}. The present paper transfers
their ideas  to the algebraic geometry over Lie algebras.

Let $A$ be a fixed Lie algebra over a field $k$. We introduce the
category of $A$-Lie algebras in Sections \ref{sec:catal} and
\ref{sec:FOL}. Sections \ref{sec:I-eofag}--\ref{sec:eqthm} are
built around the notion of a free $A$-Lie algebra $\AX$, which can
be viewed as an analogue of a polynomial algebra over a unitary
commutative ring. We introduce  a Lie-algebraic version of the
concept of an algebraic set and study connections between
algebraic sets, radical ideals of $\AX$ and coordinate algebras
(the latter can be viewed as analogues of factor-algebras of a
polynomial algebra over a commutative ring by a radical ideal).
These concepts allow us to describe the properties of algebraic
sets in two different languages:
\begin{itemize}
    \item the language of radical ideals, and
    \item the language of coordinate algebras.
\end{itemize}
One of the most important results here is Corollary
\ref{cor:as=ca}, which shows that the categories of coordinate
algebras and algebraic sets are equivalent.

In Sections \ref{sec:I-pvar}--\ref{sec:GE}, we apply some ideas of
universal algebra and model theory and introduce the notions of
$A$-prevariety, $A$-variety, $A$-quasivariety and $A$-universal
closure. We transfer some methods  of Myasnikov and Remeslennikov
\cite{alg2} to Lie algebras and solve Plotkin's problem on
geometric equivalence of Lie algebras. Our exposition is based on
a preprint by Daniyarova \cite{Dan}.

In the two final sections of this survey we describe applications
of the general theorems from Section \ref{sec:catal}--\ref{sec:GE}
to concrete classes of Lie algebras. In Section \ref{chap:II} we
survey papers \cite{AGFMLA1,AGFMLA2}:

\begin{itemize} 

\item We study the universal closure of a free metabelian Lie
algebra of finite rank $r \ge 2$ over a finite field $k$ and find
two convenient sets of axioms, $\Phi_r $ and $\Phi '_r $ for its
description; the former is written in the first order language of
Lie algebras $L$, the latter in the language $L_{\F_r}$ enriched
by constants from $\F_r$.

\item We describe the structure of finitely generated algebras
from the universal closures
     $\F_r$-$\ucl (\F_r) $ and $\ucl(\F_r)$ in languages $L_{\F_r}$ and
     $L$.

\item We prove that  in both languages $L$ and $L_{\F}$ the
universal theory of a free metabelian Lie algebra over a finite
field is decidable.
\end{itemize}

Then we apply these results to algebraic geometry over the free
metabelian Lie algebra $\F_r$, $r \ge 2$, over a finite field $k$:

\begin{itemize}
    \item We give a structural description of coordinate algebras of irreducible
    algebraic sets over  $\F_r$.
    \item We describe the structure of irreducible
    algebraic sets.
    \item We construct a  theory of dimension in the category of
    algebraic sets over  $\F_r $.
\end{itemize}

Section \ref{sec:AGOFLA} summarises the results by Daniyarova and
Remeslennikov \cite{AGOFLA} on diophantine geometry over a free Lie
algebra $F$. The objective of algebraic geometry is to classify
irreducible algebraic sets and their coordinate algebras. We believe
that the general classification problem for algebraic sets and
coordinate algebras over a free Lie algebra is very complicated,
therefore we treat only the following two cases:
\begin{itemize}
    \item algebraic sets defined by systems of equations in
    one variable;
    \item bounded algebraic sets (that is, algebraic sets contained in a finite
dimensional affine subspace of $F$, see Definition
\ref{defn:boundAS}).
\end{itemize}
In these cases,  we reduce the problem of classification of
algebraic sets and coordinate algebras to problems in  diophantine
geometry over the ground field $k$.

We refer to \cite{Bah,CK,Faith,Mal2} for background facts on Lie
algebras,  model theory, theory of categories and universal
algebra.

\tableofcontents

\section{The category of $A$-Lie algebras}
\label{sec:catal}

We work with a fixed algebra $A$ of coefficients and introduce a
notion of an $A$-Lie algebra, a Lie algebra analogue of
(associative) algebras over an associative ring.

\begin{defn} \label{defn:I-lieA}
Let $A$ be a fixed Lie algebra over a field $k$. A Lie algebra $B$
over $k$ is called an \emph{$A$-Lie algebra} if  it contains a
designated copy of $A$, which we shall usually identify with $A$.
More precisely, an $A$-Lie algebra $B$ is a Lie algebra together
with an embedding $\alpha: A \rightarrow B$. A \emph{morphism} or
\emph{$A$-homomorphism} $\varphi$ from an $A$-Lie algebra $B_1$ to
an $A$-Lie algebra $B_2$ is a homomorphism of Lie algebras which
is the identity on $A$ (or, in a more formal language,
$\alpha_1\varphi = \alpha_2$
 where $\alpha _1 $ and $\alpha _2 $ are the corresponding embeddings of
the Lie algebra $A$ into the $A$-Lie algebras $B_1$ and $B_2$).
\end{defn}

Obviously, $A$-Lie algebras and $A$-homomorphisms form a category.
In the special case $A=\left\{0\right\}$,  the category of $A$-Lie
algebras is the category of Lie algebras over $k$. Note that if
$A$ is a nonzero Lie algebra then the category of $A$-Lie algebras
does not possess a zero object.

Notice that $A$ is itself an $A$-Lie algebra.

We denote by $\hom _A (B_1, B_2)$ the set of all $A$-homomorphisms
from $B_1$ to $B_2 $,  and by $\cong _A$ the isomorphism in the
category of $A$-Lie algebras (\emph{$A$-isomorphism}). The usual
notions of free, finitely generated and finitely presented
algebras carry over to the category of $A$-Lie algebras.

We say that the set $X$ \emph{generates} an $A$-Lie algebra $B$ in
the category of $A$-Lie algebras if the algebra $B$ is generated
by the set $A \cup X$ as a Lie algebra, i. e. $B = \left< A,X
\right>$. We use notation $B =  \left< X \right>_A$.

Notice that an $A$-Lie algebra $B$ can be finitely generated in
the category of $A$-Lie algebras without  being finitely generated
as a Lie algebra.

\begin{defn} \label{defn:I-frLA}
Let $X = \left\{ x_1, \ldots, x_n \right\} $ be a finite set. An
$A$-Lie algebra $\AX = \left< x_1,\dots, x_n\right>_A$ is said to
be \emph{free} in the category of $A$-Lie algebras if, for any
$A$-Lie algebra $B=\left<b_1, \ldots, b_n \right>_A$ and any map
$\psi$ from $\AX$ to $B$ which is the identity on $A$ and
satisfies $\psi(x_i)=b_i$, $i=1, \dots, n$, the map $\psi$ extends
to an $A$-epimorphism $A[X] \longrightarrow B$. We sometimes say
that $X$ is a {\rm(}free{\rm)} base of $\AX$.
\end{defn}

A standard argument from universal algebra yields an equivalent
form of this definition:

\begin{defn} \label{defn:I-frLieAconstr}
 A \emph{free $A$-Lie algebra
with the free base $X$} is the free Lie product  of the free
{\rm(}in the category of Lie $k$-algebras{\rm)} Lie algebra $F(X)$
and algebra $A$, i.e. $\AX = A \ast F(X)$.
\end{defn}

\section{The first order language} \label{sec:FOL}

In this section we show that the category of $A$-Lie algebras is
axiomatisable.

The standard first order language $L$ of the theory of Lie
algebras over a fixed field $k$ contains a symbol for
multiplication `$\circ$', a symbol for addition `$+$', a symbol
for subtraction `$-$', a set of symbols $\left\{ k_\alpha\mid
\alpha \in k \right\}$ for multiplication  by coefficients from
the field $k$ and a constant symbol `$0$' for zero. The category
of $A$-Lie algebras requires  a bigger language $L_A$; it consists
of $L$ together with the set of constant symbols for elements in
$A$
\[
L_A  = L \cup \left\{ c_a \mid a \in A \right\}.
\]
It is clear that an $A$-Lie algebra $B$ can be treated as a model
of the language $L_A$ if the new constant symbols are interpreted
in the algebra $B$ as $c_a= \alpha(a)$.

For brevity, we sometimes omit the multiplication symbol
`$\circ$'.

Therefore the class of all $A$-Lie algebras over a field $k$ in
the language $L_A$ is given by two groups of axioms:
\begin{enumerate}
    \item The standard series of axioms that define the
    class of all Lie algebras over the field $k$.
    \item Additional axioms  \textbf{A}-\textbf{D} which describe
    the behaviour of constant symbols:
    \begin{description}
        \item[A] $0=c_0$;
        \item[B]  $c_{\alpha _1 a_1  + \alpha _2 a_2 } = k_{\alpha _1 }(c_{a_1 } )
         + k_{\alpha _2 } (c_{a_2 }) $ (for all $a_1 ,a_2  \in A$, $\alpha _1 ,\alpha _2  \in k$);
        \item[C] $c_{a_1 a_2}  = c_{a_1 } \circ c_{a_2 }$
        (for all $a_1 ,a_2  \in A$);
        \item[D] $c_a  \ne 0$ (for all nonzero $a \in A$).
     \end{description}
\end{enumerate}

Axioms \textbf{A}, \textbf{B}, \textbf{C} and \textbf{D} imply
that in an $A$-Lie algebra $B$ constant symbols $c_a$ are
interpreted as distinct elements and form a subalgebra of $B$
isomorphic to $A$.

\section{Elements of algebraic geometry}
\label{sec:I-eofag}

Our next objective is to introduce Lie algebraic counterparts to
the classical concepts of algebraic geometry.

Let $A$ be a fixed Lie algebra over a field $k$.
Let $X =\left\{ x_1,\ldots ,x_n \right\} $ be a finite set of
variables and $\AX$ be the free $A$-Lie algebra with the base $X$.
We view $\AX$ as an analogue of a polynomial algebra in finitely
many variables over a unitary commutative ring. We think of
elements of $\AX$ as polynomials with coefficients in $A$. We use
functional notation
\[
f = f(x_1 ,\ldots ,x_n ) = f(x_1 ,\ldots ,x_n ,a_1 ,\ldots ,a_r )
\]
thereby expressing the fact that the Lie polynomial $f$ in $\AX$
involves {\it variables} $x_1,\dots,x_n$ and, if needed, {\em
constants} $a_1,\dots,a_r \in A$. A formal equality $f=0$   can be
treated, in an obvious way, as an equation over $A$. Therefore
every subset $S \subset \AX$ can be treated as a \emph{system of
equations} with coefficients in $A$. In parallel with  the
commutative case, the set of solutions of $S$ depends on the
algebra used to solve the system. We are specially interested in
\emph{diophantine problems}, that is, solving systems in $A$,  but
for the time being we work in a more general situation and solve
$S$ in an arbitrary $A$-Lie algebra $B$.


Let $B$ be an $A$-Lie algebra and  $S$  a subset of $\AX$. Then
the set $B^n =  \left\{ (b_1 ,\ldots ,b_n ) | \;b_i  \in B
\right\}$ is called the \emph{affine $n$-dimensional space over
the algebra $B$}.

A point  $p = (b_1,\ldots , b_n ) \in B^n $  is called a
\emph{root of a polynomial} $f \in \AX$ if
\[
 f(p) = f(b_1 ,\ldots ,b_n, a_1 ,\ldots,a_r ) = 0.
\]
We also say that the polynomial $f$ \emph{vanishes} at the point
$p$. A point $p \in B^n $ is called  \emph{a root or a solution of
the system} $S \subseteq \AX$  if every polynomial from $S$
vanishes at $p$.

\begin{defn}
Let $B$ be an $A$-Lie algebra and let $S$ be a subset
of $\AX$. Then the set
\[
V_B (S) = \left\{ p \in B^n  | \; f(p) = 0 \ \ \forall \, f \in S
\right\}
\]
is called the {\rm(}\emph{affine}{\rm)} \emph{algebraic set over
$B$ defined by $S$}.
\end{defn}

\begin{defn}
Let $B$ be an $A$-Lie algebra, $S_1$ and $S_2$  subsets of $\AX$.
Then the systems $S_1$ and $S_2$ are called \emph{equivalent over
$B$} if \/ $V_B (S_1 ) = V_B (S_2)$. A system $S$ is called
\emph{inconsistent over $B$} if $V_B (S) = \emptyset$ and
\emph{consistent} otherwise.
\end{defn}

{\rm
\begin{expl}[Typical examples of algebraic sets] \label{expl:1} \
\begin{enumerate}

\item Every element $a \in A$ forms an algebraic set, $\left\{ a
\right\}$. Indeed if you take $S = \left\{ x - a \right\}$ then
         $V_B (S) = \left\{ a \right\}$. In this example $n = 1$
         and $X = \left\{ x \right\}$.
\item Every element $\left\{ (a_1, \ldots, a_n) \right\} \in A^n$
is an algebraic set: if you take $S = \left\{ x_1 - a_1 , \dots,
x_n-a_n \right\}$ then  $V_B (S) = \left\{ (a_1, \ldots, a_n)
\right\}$.
    \item The centraliser $C_B (M)$ of an arbitrary set of elements $M$ from
         $A$ is an algebraic set defined by the system $S = \left\{ x
         \circ m \mid m \in M \right\}$.
    \item The whole affine space $B^n$ is the algebraic set defined by the system $S = \left\{
         0\right\}$.
    \item Let $A=\left\{0 \right\}$. Then the empty set
         $\emptyset$  is not algebraic, since every algebraic set in
         $B^n$ contains the point $(0,\ldots,0)$.
\end{enumerate}
\end{expl}
}

We need more definitions.

A polynomial $f \in \AX$ is called a \emph{consequence of the
system} $S \subseteq \AX$ if $ V(f) \supseteq V(S)$.

Let $Y$ be an arbitrary  {\rm(}not necessarily algebraic{\rm)}
subset of $B^n $. The set $$ \rad_B (Y) =  \left\{ f \in \AX |\;
f(p) = 0 \ \ \forall \, p \in Y \right\}$$ is called the
\emph{radical} of the set $Y$. If \/ $Y=\emptyset$ then, by the
definition, its radical is the algebra $\AX$.

If $Y$ is an algebraic set ($Y = V_B (S)$) then we also refer to
its radical as the radical of the system of equations $S$: $
\rad_B (S) = \rad_B (V_B (S))$, i.e.
 {\rm
\begin{itemize}
    \item a polynomial $f \in \AX$ is a consequence of a system
$S$ if and only if $f \in \rad_B (S)$;
    \item a polynomial $f$ is
a consequence of a system $S$ if and only if the system $S' = S
\cup \left\{ f \right\}$ is equivalent to $S$.
\end{itemize}}
Therefore, $\rad_B (S)$ is the maximal (by inclusion) system of
equations equivalent to $S$.

\begin{prop} \label{prop:I-rad=cons}
 The radical of a set is an ideal of the algebra $\AX$.
\end{prop}
\begin{proof} Let $f,g \in \rad_B (Y)$ and $h \in \AX$ and let $y \in Y$.
By  definition, $$(\alpha f+\beta g)(y)=\alpha f(y)+\beta g(y)=0$$
and $(hf)(y)=h(y)\cdot f(y)=h(y)\cdot 0=0$, where $\alpha, \beta
\in k$.
\end{proof}

\begin{lem} \label{lem:rad}\
\begin{enumerate}
    \item The radical of a system $S \subseteq \AX$ contains the
            ideal \/ $\verb"id" \left< S \right>$ generated by the set $S$,
            {\rm ${\rm \rad}_B (S) \supseteq \verb"id" \left< S \right>.$}
    \item Let\ $Y_1$ and $Y_2$  be  subsets of $B^n$ and $S_1$, $S_2$   subsets of $\AX$.

If\ $Y_1  \subseteq Y_2$  then {\rm ${\rm \rad_B} (Y_1) \supseteq
{\rm \rad_B (Y_2)}$}.

If\ $S_1  \subseteq S_2$ then {\rm ${\rm \rad_B}(S_1) \supseteq {\rm
\rad_B (S_2)}$}.

    \item For any family of sets $\left\{\, Y_i | \; i \in I
\,\right\}$, $Y_i \subseteq B^n $ we have {\rm
\[
{\rm \rad}_B \left(\bigcup\limits_{i \in I} {Y_i }\right) =
\bigcap\limits_{i \in I} {{\rm \rad}_B (Y_i )}.
\]}

    \item An ideal $I$ of the algebra $\AX$ is the radical of an
algebraic  set over $B$ if and only if \/ {\rm $ {\rm \rad_B} (V_B
(I)) = I$}.
    \item A set $Y \subseteq B^n$ is algebraic over $B$ if
    and only if \/ {\rm $ V_B ({\rm \rad} _V (Y)) = Y. $}
    \item Let\/ $Y_1 ,Y_2  \subseteq B^n $ be two algebraic sets, then
{\rm
\[
Y_1  = Y_2 \hbox{ \textit{if and only if} }  \rad_B (Y_1 ) =
\rad_B (Y_2 ).
\]
}
\end{enumerate}
Therefore the radical of an algebraic set describes it uniquely.
\end{lem}

\begin{proof} The proofs immediately follow from the definitions.
As an example, we prove the fourth statement. If \/$\rad_B (V_B
(I)) = I$ then $I$ is obviously a radical. \emph{Vice versa} if
$I$ is the radical of an algebraic set then there exists a system
$S$ such that $I=\rad_B(S)$. Then $V_B(I)=V_B(\rad_B(S))=V_B(S)$,
consequently, $\rad_B(V_B(I))=\rad_B(V_B(S))=\rad_B(S)=I$.
\end{proof}

Another crucial concept notion of algebraic geometry is that of
\emph{coordinate algebra}.

\begin{defn} \label{defn:I-cooralg}
Let $B$ be an $A$-Lie algebra, $S$  a subset of $\AX$ and
$Y\subseteq B^{n}$  the algebraic set defined by the system $S$.
Then the factor-algebra\\[.3ex]

\centerline{{\rm $ \Gamma _B (Y) = \Gamma _B (S) =
{\raise0.7ex\hbox{${\AX}$} \!\mathord{\left/
 {\vphantom {{\AX} {\rad_B (Y)}}}\right.\kern-\nulldelimiterspace}
\!\lower0.7ex\hbox{${\rad_B (Y)}$}} $}}

\medskip \noindent is called the \emph{coordinate algebra} of the
algebraic set $Y$ {\rm(}or of the system $S${\rm)}.
\end{defn}

\begin{rem} \label{rem:I-catcoordalg}
Observe that if a system  $S$ is inconsistent over $B$ then
$\Gamma _B (S) = 0$. Notice that the  coordinate algebras of
consistent  systems of equations are $A$-Lie algebras and form a
full subcategory of the category of all $A$-Lie algebras.
\end{rem}

The coordinate algebra of an algebraic set can be viewed from a
different perspective, as an algebra of polynomial functions.
Indeed if $Y = V_B (S)$ is an algebraic set in $B^n$ then the
coordinate algebra $\Gamma _B (Y)$ can be identified with the
$A$-Lie algebra of all polynomial functions on $Y$; the latter are
the functions from $Y$ into $B$ of the form
\[
\bar{f}: Y \rightarrow B, \quad p \rightarrow f(p)\;\;\; (p\in Y),
\] where $f \in A\left[x_1,\ldots, x_n\right]$ is a polynomial.

It is clear that two polynomials $f,g \in \AX$ define the same
polynomial function if and only if $f - g \in \rad_B (Y)$. The set
of all polynomial functions $P_B (Y)$ from  $Y$ to $B$ admits the
natural structure of an $A$-Lie algebra. We formulate our
observation as the following proposition.

\begin{prop} \label{prop:I-coralg-pol}
Let\/ $B$ be an $A$-Lie algebra and\/ $Y$ a nonempty algebraic
set. Then the coordinate algebra $\Gamma _B (Y)$ of\/ $Y$ is the
$A$-Lie algebra of all polynomial functions on\/ $Y$: $$ \Gamma _B
(Y) \cong _A P_B (Y).$$
\end{prop}

\begin{expl}
If $a \in A$ and $Y = \left\{ a \right\}$ then $\Gamma _B (Y)
\cong A$.
\end{expl}

Similarly to the commutative case, points of an algebraic set can
be viewed as certain Lie algebra homomorphisms.

\begin{prop} \label{prop:algsashom}
Every algebraic set\/ $Y$ over $B$ can be identified with the set
{\rm $\hom_A (\Gamma _B (Y),B)$} {\rm(}see Section {\rm
\ref{sec:catal}} for notation{\rm)} by the rule {\rm
$$
\theta: \hom_A
(\Gamma _B (Y),B) \leftrightarrow Y.
$$}
Consequently, for any
algebraic set $Y$ there is a one-to-one correspondence between the
points of $Y$  and $A$-homomorphisms from $\Gamma _B(Y)$ to $B$.
\end{prop}
\begin{proof} Indeed the coordinate algebra $\Gamma _B (Y)$ is the
factor-algebra $$ \Gamma _B (Y) = {\raise0.7ex\hbox{${\AX}$}
\!\mathord{\left/
 {\vphantom {{\AX} {\rad_B (Y)}}}\right.\kern-\nulldelimiterspace}
\!\lower0.7ex\hbox{${\rad_B (Y)}$}}$$ for $X = \left\{ x_1, \ldots
,x_n \right\}$. Defining an $A$-homomorphism $\varphi \in \hom_A
(\Gamma _B (Y),B)$ is equivalent to defining the images in $B$ of
the elements of $X$, that is, to fixing a point $(b_1 ,\dots ,b_n)
\in B^n $. Therefore we set $\theta (\varphi)=(b_1 ,\dots ,b_n)$,
where $(b_1 ,\ldots ,b_n)$ is the image of $X$ under $\varphi$.
The point $(b_1 ,\ldots ,b_n)$ must satisfy the  condition:
\[
f(b_1 ,\ldots,b_n ) = 0 \mbox{ for all }   f \in \rad_B (Y).
\]
Clearly, this condition holds only at the points of the algebraic
set $Y$. Obviously, distinct homomorphisms correspond to distinct
points and every point in $Y$ has a pre-image in $\hom_A (\Gamma
_B (Y),B)$.
\end{proof}

Developing these  ideas further, we come to the following method
for computation of the radical $\rad_B (Y)$ of an algebraic set
$Y$.

Let $p \in B^n$ be an arbitrary point and denote by $\varphi _p $
the  $A$-homomorphism
\begin{equation} \label{eq:I-observ}
 \varphi _p: \AX \rightarrow B,\ \ f \in \AX,\ \  \varphi _p (f) = f(p) \in B.
\end{equation}

Clearly, $\rad_B(\left\{p\right\}) = \ker \varphi _p $. In view of
Lemma \ref{lem:rad} we have $$ \rad_B (Y) = \mathop \bigcap
\limits_{p \in Y} \ker \varphi _p.$$

 This equality clarifies
the structure of the radical of an algebraic set and the structure
of its coordinate algebra. Indeed, by Remak's Theorem we have an
embedding $$ \Gamma _B (Y) \rightarrow \prod\limits_{p \in Y}
{{\raise0.7ex\hbox{${\AX}$} \!\mathord{\left/
 {\vphantom {{\AX} {\ker \varphi _p }}}\right.\kern-\nulldelimiterspace}
\!\lower0.7ex\hbox{${\ker \varphi _p }$}}}. $$ The factor-algebra
$ {\raise0.7ex\hbox{${\AX}$} \!\mathord{\left/
 {\vphantom {{\AX} {\ker \varphi _p }}}\right.\kern-\nulldelimiterspace}
\!\lower0.7ex\hbox{${\ker \varphi _p }$}}$ is isomorphic to $\im\,
\varphi _p$ and thus imbeds into $B$. This implies that the
coordinate algebra $\Gamma _B (Y)$ imbeds into a cartesian power
of the algebra $B$; we state this observation as the following
proposition.

\begin{prop} \label{prop:I-coormapstoB}
The coordinate algebra of an algebraic set over $B$  imbeds into
an unrestricted cartesian power of $B$, $$ \Gamma _B (Y)
\hookrightarrow B^Y.$$
\end{prop}

In particular, Proposition \ref{prop:I-coormapstoB} implies that
all the identities and quasi-identities which are true in $B$ are
also true in $\Gamma _B (S)$ (see Section \ref{sec:I-uncl} for
definitions). In particular, if $B$ is either an abelian, or
metabelian, or nilpotent Lie algebra, the coordinate algebra
$\Gamma_B(Y)$ of an arbitrary algebraic set $Y$ over $B$ is
abelian, or metabelian, or nilpotent, respectively; in the latter
case, the nilpotency class of $\Gamma_B(Y)$ does not exceed the
nilpotency class of $B$.

\section{The Zariski topology} \label{sec:I-zartop}

To introduce the Zariski toppology on the $n$-dimensional affine
space $B^n $, we formulate first the following auxiliary result.

\begin{lem} \label{lem:zar}
If  $A$ is a nonzero Lie algebra then
\begin{enumerate}
    \item the empty set\/ $\emptyset $ is an algebraic set over any nonzero $A$-Lie
algebra $B$: it suffices to take  $S = \left\{ a \right\}$ for
some non-zero $a \in A$;
    \item the whole affine space $B^n$ is an algebraic
set {\rm(}see Example {\rm \ref{expl:1})};
    \item the intersection of a family
of algebraic sets over $B$ is also an algebraic set over $B$:
$$
\bigcap\limits_{i \in I} {V_B (S_i )}  = V_B
\left(\bigcup\limits_{i \in I} {S_i }\right) \ \mbox{ for }  \ S_i
\in \AX.
$$
\end{enumerate}
\end{lem}

In view of Lemma \ref{lem:zar}, we can define a topology in $B^n$
by taking algebraic sets in $B^n$ as a sub-basis for the
collection of closed sets. We call this topology the {\it Zariski
topology}. For later use we introduce the following notation for
three families of subsets of $B^n$:
\begin{itemize}
    \item $\Upsilon$ is the collection of all algebraic sets
over $B$ (the sub-basis of the topology).
    \item  $\Upsilon_1$ is the collection of all
finite unions of  sets from  $\Upsilon$ (the basis of the
topology).
    \item $\Upsilon_2$ is the collection of all intersections  of sets
from $\Upsilon_1$, i.e.\ $\Upsilon_2$ is the set of all Zariski
closed subsets of the space $B^n $.
\end{itemize}

\begin{defn} \label{defn:I-noeth}
An $A$-Lie algebra $B$ is called \emph{$A$-equationally
Noetherian} if for any $n \in \mathbb{N}$ and any system
 $S \subseteq A\left[x_1 ,\ldots,x_n \right]$
 there exists a finite subsystem  $S_0  \subseteq S$
such that $V_B (S) = V_B (S_0 )$.
\end{defn}

Recall that a topological space $(T,\tau)$ is called
\emph{Noetherian} if and only if every strictly descending chain
(with respect to the inclusion) of closed subsets terminates.
Provided that a sub-basis $\sigma$ of $\tau$ is closed under
intersections, one can give an equivalent formulation:  a
topological space $(T,\tau)$ is  Noetherian if and only if every
strictly descending chain  of subsets from $\sigma$ terminates
\cite{Top}.

\begin{lem}
The algebra $B$ is $A$-equationally Noetherian if and only if, for
every positive integer $n$, the affine space $B^n$ is Noetherian.
\end{lem}

\begin{proof} Let us assume for the time being  that $B$ is $A$-equationally
Noetherian. Consider a descending chain of algebraic sets: $Y_1
\supset Y_2 \supset \ldots \supset Y_t \supset \ldots$. By Lemma
\ref{lem:rad}, the radicals of these sets form an increasing
chain: $ \rad_B (Y_1 ) \subset \rad_B (Y_2 ) \subset \dots$. Let
$$S = \bigcup \limits_i {\rad_B (Y_i )}.$$ Since the system $S$ is
equivalent to its finite subsystem $S_0$ the chain of radicals
contains only finite number of pairwise distinct ideals, therefore
the both chains terminate after finally many steps.

Now suppose that, for every positive integer $n$, the space $B^n $
is topologically Noetherian. We wish to show that the $A$-Lie
algebra $B$ is $A$-equationally Noetherian. Let  $S \subseteq
A\left[x_1 ,\dots,x_n \right]$ be an arbitrary system of equations
and  $s_1$  a polynomial from $S$. If the system $S$ is equivalent
(over $B$) to its subsystem $S_0 = \left\{ s_1 \right\}$, then the
statement is straightforward. Otherwise, there exists an element
$s_2 \in S \smallsetminus \left\{ s_1 \right\}$ such that $V_B
(\left\{ s_1 \right\}) \supset V_B (\left\{ s_1 ,s_2 \right\} )$;
we continue the construction in an obvious way and build a
decreasing chain
\[
V_B (\left\{ s_1 \right\}) \supset V_B (\left\{ s_1 ,s_2
\right\})\supset \ldots \supset V_B (\left\{ s_1 , \ldots,
s_r)\right\})\supset \ldots
\]
Since the space $B^n $ is topologically Noetherian the chain
contains only finite number of pairwise distinct sets. Therefore
the system $S$ is equivalent (over $B$) to a finite subsystem.
\end{proof}


A closed set $Y$ is called \emph{irreducible} if \/ $Y=Y_1 \cup
Y_2$, where $Y_1$ and $Y_2$ are closed, implies that either
$Y=Y_1$ or $Y=Y_2$.

\begin{thm} \label{thm:1-1}
Any closed subset $Y$ of $B^n$ over an $A$-equationally Noetherian
$A$-Lie algebra $B$ can be expressed as a finite union of
irreducible algebraic sets: $$Y = Y_1 \cup \cdots \cup Y_l.$$ This
decomposition is unique if we assume, in addition, that $Y_i
\nsubseteq Y_j$ for $i\ne j$; in that case, $Y_i$  are referred to
as the irreducible components of $Y$.
\end{thm}

Proof is standard, see, for example, \cite{RH}.

The dimension of an irreducible algebraic set $Y$ is defined in
the usual way.

Let $Y$ be an irreducible algebraic set. A supremum, if exists, of
all integers $m$ such that there exists a chain of irreducible
algebraic sets
\[
Y = Y_0  \varsupsetneq Y_1  \varsupsetneq \ldots \varsupsetneq Y_m
\]
is called the \emph{dimension} of $Y$ and is denoted by $\dim
(Y)$. If the supremum does not exist then, by definition, we set
$\dim (Y)=\infty$.

We can generalise this definition when $Y$ is an arbitrary (not
necessarily irreducible) algebraic set over an $A$-equationally
Noetherian $A$-Lie algebra $B$: we define its dimension $\dim (Y)$
as the supremum of the dimensions of its irreducible components.

\section{$A$-Domains}

For groups, the notion of a \emph{domain} (a group that has no
zero divisors) was introduced in \cite{alg1}; it is useful for
formulation of irreducibility criteria for algebraic sets over
groups. Here, we introduce a similar concept for Lie algebras.

\begin{defn}
Let $B$ be an $A$-Lie algebra, $x \in B$.
 \bi
\item The principal ideal of the subalgebra $\left< A,x \right>$
generated by $x$ is called the \emph{$A$-relative ideal generated
by} $x$ and denoted by $\left< x \right>^A$.

 \label{defn:I-0div}

\item A nonzero element $x \in B$ is called an \emph{$A$-zero
divisor} if there exists $y \in B$, $y \ne 0$ such that $\left< x
\right>^A \circ \left< y \right>^A = 0$.

\item  \label{defn:I-dom} An $A$-Lie algebra $B$ is  an
\emph{$A$-domain} if it contains no nonzero $A$-zero divisors. \ei
\end{defn}

\begin{expl} \label{ex:1}
Let $A$ be a nilpotent Lie algebra. Then every element $x \in A$
is a zero divisor.

Indeed, take a nonzero element $y$ from the center of $A$. Then
$\verb"id"\left< y \right>$ is a $k$-vector space with the basis
$\left\{ y\right\}$. The pair $(x,y)$ is a pair of zero divisors.

Recall that the \emph{Fitting's radical} of a Lie algebra $A$ is
the ideal generated by the set of all elements from the nilpotent
ideals of $A$. Obviously, if
 \label{ex:2}
$A$ is an arbitrary Lie algebra over $k$, then every nonzero
element $x$ of the Fitting's radical \/ {\rm $\Fit (A)$}  of the
algebra $A$ is a zero divisor since  the ideal $\verb"id"\left< x
\right>$ is  nilpotent {\rm(}see {\rm\cite{AGFMLA1})}.
\end{expl}

The next lemma shows that in an $A$-domain $B$ every closed subset
of the space $B^n$ is algebraic.

\begin{lem} \label{lem:clsdo}
Let $B$ be an $A$-domain. Then $\Upsilon=\Upsilon_1=\Upsilon_2$.
\end{lem}

\begin{proof} By Lemma \ref{lem:zar} it will suffice to show
that if $B$ is an $A$-domain then any finite union of algebraic
sets is also algebraic.

Let $ S_1 ,S_2 \subset \AX$ be two consistent systems of equations
and assume that $V_B (S_1)$ and $V_B (S_2 )$ are the algebraic
sets defined by the systems $S_1$ and $S_2$, respectively. Then
$V_B (S_1) \cup V_B (S_2 ) = V_B (S)$, where here
\[
S = \left\{ f_1  \circ f_2 |\;f_i  \in \left< s_i \right>^A,\;s_i
\in S_i, \;i = 1,2 \right\}.
\]
\end{proof}

\begin{lem} Let $B$ be an
$A$-domain and $Y$ be an arbitrary subset of $B^n $. Then the
closure of \/ $Y$ in the Zariski topology coincides with\/ {\rm
$V_B ( \rad_B (Y))$}.
\end{lem}
\begin{proof} Clearly the set $V_B (\rad_B (Y))$ is closed and contains
$Y$. We show that $V_B (\rad_B (Y))$ is contained in every closed
set $Z$ such that $Y \subseteq Z$. According to Lemma
\ref{lem:rad}, $\rad_B (Y) \supseteq \rad_B (Z)$ and thus $V_B
(\rad_B (Y)) \subseteq V_B (\rad_B (Z))$. By Lemma
\ref{lem:clsdo}, every closed set in $B^{n}$ is algebraic over
$B$, hence $V_B (\rad_B (Z)) = Z$ and the statement follows.
\end{proof}

\section{The category of algebraic sets} \label{sec:I-AS}

In this section we introduce the category $\AS_{A,B}$ of algebraic
sets over an $A$-Lie algebra $B$. Throughout this section we
assume that $B$  is  an $A$-Lie algebra,  $\AX = A\left[x_1,
\ldots,x_n \right]$ and that $B^n$ is the affine $n$-space over
$B$.

\begin{defn} \label{defn:I-polmap}
Objects of $\AS_{A,B}$ are algebraic sets in all affine $n$-spaces
$B^n$, $n\in \mathbb{N}$. If $Y \subseteq B^n$ and $Z \subseteq
B^m$ are algebraic sets, then a map $ \psi: Y \rightarrow Z$ is a
\emph{morphism} if there exist $f_1,\dots,f_m \in
A\left[x_1,\ldots ,x_n \right]$
 such that, for any $(b_1,\ldots,b_n) \in Y$,
\[
\psi(b_1,\dots,b_n)=(f_1(b_1,\dots,b_n),\ldots,f_m(b_1,\dots,b_n))
\in Z.
\]
Occasionally we  refer to morphisms as \emph{polynomial maps}.
 \label{defn:I-AlgSet}

We denote by $\hom(Y,Z)$ the set of all morphisms from $Y$ to $Z$.

 \label{defn:I-asisom} Following usual conventions of category
theory, algebraic sets $Y$ and $Z$ are called \emph{isomorphic} if
there exist morphisms
\[
\psi: Y \rightarrow Z \hbox{ and } \theta: Z \rightarrow Y
\]
such that \[ \theta \psi  = \verb"id"_Y \hbox{ and } \psi \theta =
\verb"id"_Z.
\]
\end{defn}

\begin{lem}
Let $Y \subseteq B^n$ and $Z \subseteq B^m$ be  algebraic sets
over $B$ and $\psi$  a morphism from $Y$ to $Z$. Then
\begin{enumerate}
    \item  $\psi$ is a continuous map in the
Zariski topology;
    \item if\/ $Y$ is an irreducible  algebraic set and $\psi$ is an epimorphism then $Z$ is also irreducible.
    In particular, the {\rm(}ir{\rm)}reducibility of algebraic sets is preserved by isomorphisms.
\end{enumerate}
\end{lem}
\begin{proof}
(1) By  definition, $\psi$ is continuous if and only if the
pre-image $Y_1 \subseteq B^n$ of any algebraic set $Z_1 \subseteq
B^m$  is also algebraic. Let $Z_1 = V_B (S)$ then
\[
\psi ^{ - 1} (Z_1 ) =  \left\{ p \in B^n  |\;S(f_1 (p), \ldots,f_m
(p)) = 0 \right\}.
\]
Clearly this subset is algebraic over $B$.

(2) Suppose that  $Z$ is reducible and $Z = Z_1  \cup Z_2 $. Then
$Y = \psi^{-1}(Z_1) \cup \psi^{-1}(Z_2)$ and $\psi^{- 1} (Z_1 )$,
$\psi ^{ - 1} (Z_2 )$ are proper closed subsets of $Y$.
\end{proof}

\begin{lem} \label{lem:I-dir-alg}
Let $Z_1 \subseteq B^n $ and $Z_2  \subseteq B^m $ be algebraic
sets over $B$. Then the set $Z_1  \times Z_2  \subseteq B^{n + m}
$ is also algebraic over $B$.
\end{lem}
\begin{proof} Let $Z_1  = V_B (S_1 )$, $S_1 \in \AX$, $X = \left\{ x_1,
\ldots,x_n \right\}$ and $Z_2  = V_B (S_2)$, $S_2  \in
A\left[Y\right]$, $Y = \left\{ y_1 ,\ldots,y_m \right\}$ and let
$X$ and $Y$ be disjoint. Observe that the set $Z_1 \times Z_2 $ is
defined by the system of equations $S_1 \cup S_2 \subseteq
A\left[X \cup Y\right]$. \end{proof}

\begin{expl} \label{expl:I-dim1}
Let $A$ be a Lie $k$-algebra with  trivial multiplication and
assume that $B = A$. An elementary re-interpretation of some basic
results of linear algebra  yields the following results.
\begin{enumerate}
    \item Every consistent system of equations over $A$ is equivalent
to a triangular system of equations.
    \item The morphisms in the category {\rm $\AS_{A,A}$} are affine transformations.
    \item Every algebraic set $Y \subseteq A^n$ is
isomorphic to an algebraic set of the form $ (\underbrace
{A,A,\ldots,A}_s,0,\ldots,0), \ \ 0 \le s \le n$.
    \item  Every coordinate algebra $\Gamma (Y)$ is  $A$-isomorphic to
$$A \oplus \verb"lin"_k \left\{ x_1 ,\ldots,x_s \right\},$$ where $0
\le s \le n$, and \/ $\verb"lin"_k \left\{ x_1 ,\ldots,x_s
\right\}$ is the linear span of the elements $\left\{ x_1
,\ldots,x_s \right\}$ over $k$.
\end{enumerate}
\end{expl}

\section{The Equivalence Theorem} \label{sec:eqthm}

One of the main problems in algebraic geometry is the problem of
classification of algebraic sets up to isomorphism. In this
section we prove that this problem  is equivalent to
classification of coordinate algebras.

Denote by  $\CA_{A,B}$  the category of all coordinate algebras of
algebraic sets from $\AS_{A,B}$ (morphisms in $\CA_{A,B}$ are
$A$-homomorphisms). As we observed in Section \ref{sec:I-eofag},
coordinate algebras of nonempty algebraic sets over $B$ form a
full subcategory of the category of $A$-Lie algebras (see Remark
\ref{rem:I-catcoordalg}).

Notice that if the empty set is an algebraic set then the zero
algebra is a coordinate algebra, and {\rm
$\hom_A(0,\Gamma_B(S))=\hom_A(\Gamma_B(S),0)=\emptyset$}.

 For the time being we  declare (and later show) that the
categories $\AS_{A,B}$ and $\CA_{A,B}$ are equivalent but not
isomorphic.

A pair $(X,S)$, where here $X = \left\{ x_1 ,\dots,x_n \right\}$
and $S \subseteq \AX$, where $S$ is a radical ideal, is called a
\emph{co-presentation} of the coordinate algebra $$ \Gamma _B (S)
= {\raise0.7ex\hbox{${\AX}$} \!\mathord{\left/
 {\vphantom {{\AX} {\rad_B (S)}}}\right.\kern-\nulldelimiterspace}
\!\lower0.7ex\hbox{${\rad_B (S)}$}}.$$

Let $(X,S_1 )$, $(Y,S_2 )$ be co-presentations, $X = \left\{ x_1
,\ldots,x_n \right\}$, $S_1  \subseteq \AX$, $Y= \left\{ y_1
,\ldots,y_m \right\}$, $S_2 \in A\left[Y\right]$. An
$A$-homomorphism $\varphi: \AX \rightarrow A\left[Y\right]$ so
that {\rm $\varphi (\rad_B (S_1 )) \subseteq \rad_B (S_2 )$} is
called a morphism from the co-presentation $(X,S_1 )$ to the
co-presentation $(Y,S_2 )$.

Naturally, co-presentations $(X,S_1 )$ and $(Y,S_2 )$ are called
\emph{isomorphic} if the respective coordinate algebras $\Gamma _B
(S_1)$ and $\Gamma _B (S_2 )$ are $A$-isomorphic.

The collection of all co-presentations together with
 the morphisms defined above form a category, which we call
the category of co-presentations of coordinate algebras and denote
it by \/ {\rm $\RCA_{A,B}$}.

\begin{thm} \label{thm:I-equiv}
The categories {\rm $\AS_{A,B}$} and \/ {\rm $\RCA_{A,B}$} are
isomorphic.
\end{thm}
\begin{proof} We construct
two contravariant functors
\[
 \Ff: \AS_{A,B}  \rightarrow  \, \RCA_{A,B} \hbox{ and } \Gg: \RCA_{A,B}  \rightarrow \AS_{A,B}
 \]
such that $\Ff\Gg  = \id _{\RCA_{A,B}}$  and  $\Gg\Ff =
\id_{\AS_{A,B}}$.

We define the functors $\Ff$ and $\Gg$ first on the objects and
then on the morphisms of respective categories.

 Every algebraic set $Y$, as well as
every co\/-presentation, is defined by the cardinality $n$ of the
set of variables $X = \left\{ x_1 ,\ldots,x_n \right\}$ and by the
radical $S=\rad(Y)\subset\AX$ treated as a system of equations.
We, therefore, set:
\[
\Ff(V_B(S)) = (X,S),\quad \Gg((X,S)) = V_B (S).
\]

Next we define the functors on the morphisms. To this end, let
$Z_1  = V_B (S_1 )$ and $Z_2  = V_B (S_2 )$ be  algebraic sets
over $B$, where  $S_1 \subseteq \AX$, $X = \left\{ x_1,\dots, x_n
\right\}$ and $S_2 \in A \left[ Y \right]$, $Y = \left\{ y_1
,\ldots, y_m \right\}$ are radicals. Suppose that $(X,S_1 )$,
$(Y,S_2 )$ are the respective co-presentations. By  definition, a
contravariant functor is a morphism-reversing functor, i.e. {\rm
\begin{gather} \notag
\begin{split}
  \psi  \in \hom(Z_1 ,Z_2 ) & \rightarrow    \Ff(\psi ) \in \hom((Y,S_2 ),(X,S_1 )) \\
    \varphi  \in \hom((Y,S_2 ),(X,S_1 ))  & \rightarrow  \Gg(\varphi)\in \hom(Z_1 ,Z_2 ).
\end{split}
\end{gather}
} Choose $\psi  \in \hom(Z_1 ,Z_2)$ and polynomials $f_1, \ldots,
f_m  \in A\left[x_1 ,\ldots,x_n \right]$ so that
\[
\psi (b_1 ,\ldots,b_n ) = (f_1 (b_1 ,\ldots,b_n ),\ldots,f_m (b_1,
\ldots, b_n)),\quad (b_1 ,\ldots,b_n ) \in Z_1.
\]
We define $\Ff(\psi): A\left[Y\right] \rightarrow \AX$ by
$\Ff(\psi)(y_i) = f_i (x_1 ,\ldots,x_n )$, $i = 1,\ldots,m$. By
definition, $\Ff(\psi) (S_2 ) \subseteq S_1$ and therefore
$\Ff(\psi)$ is a morphism in the category $\RCA_{A,B}$.

Now suppose that $\varphi  \in \hom((Y,S_2 ),(X,S_1 ))$, i.e.
$\varphi: A\left[Y\right] \rightarrow \AX$ is an $A$-homomorphism
satisfying $\varphi (S_2 ) \subseteq S_1$. We define the image
$\Gg(\varphi ): Z_1 \rightarrow Z_2$ of $\varphi$ under $\Gg$ by
polynomial maps $f_1, \ldots, f_m \in A\left[x_1, \ldots, x_n
\right]$, where $f_1, \ldots, f_m$ are the images of elements of
$Y$ under $\varphi$. It is easy to check the inclusion
$\Gg(\varphi)(Z_1 ) \subseteq Z_2$; the equalities
\[
\Ff\Gg = \verb"id" _{\RCA_{A,B} } \hbox{ and } \Gg\Ff =
\verb"id"_{\AS_{A,B}}
\]
follow from definition.

To complete the proof of the theorem we need to show that $\Ff$
and $\Gg$ are in fact functors. By definition, a contravariant
functor $\Ff$  satisfies two following conditions:
\begin{itemize}
    \item for an arbitrary algebraic
    set $Z$: $\Ff(\verb"id"_Z ) = \verb"id"_{\Ff(Z)}$  and
    \item for any two morphisms $\psi $ and $\theta$ in $\AS_{A,B}$:
    $\Ff(\psi \theta ) = \Ff(\theta )\Ff(\psi )$.
\end{itemize}
The verification of these conditions for $\Ff$ and $\Gg$ is
straightforward and left to the reader.
\end{proof}

\begin{cor}
Two algebraic sets are isomorphic if and only if the respective
co-presentations are. Two co-presentations are isomorphic if and
only if the respective coordinate algebras are $A$-isomorphic.
\end{cor}

In informal terms, the categories $\RCA_{A,B}$ and $\CA_{A,B}$
look very much alike. The correspondence between the objects and
the morphisms of the categories $\RCA_{A,B}$ and $\CA_{A,B}$,
therefore, establishes a correspondence between the categories
$\AS_{A,B}$ and $\CA_{A,B}$; however, this correspondence  is not
one-to-one: the same coordinate algebra corresponds to different
algebraic sets because it has  different co-presentations.

\begin{cor} \label{cor:as=ca}
The category \/ {\rm $\AS_{A,B}$} of algebraic sets over an
$A$-Lie algebra $B$ is equivalent to the category \/ {\rm
$\CA_{A,B}$} of coordinate algebras over $B$.
\end{cor}

We note also a corollary from the proof of Theorem {\rm
\ref{thm:I-equiv}}:

\begin{cor} \label{cor:I-epi}
Let $Y$ and $Z$ be algebraic sets over an algebra $A$ and
\/$\Gamma(Y)$ and \/$\Gamma(Z)$ the respective coordinate
algebras. Then we have a one-to-one correspondence between \/{\rm
$\hom (Y,Z)$} and \/ {\rm $\hom (\Gamma(Y),\Gamma(Z))$}. Moreover,
every embedding of algebraic sets $Y \subseteq Z$ corresponds to
an $A$-epimorphism $\varphi:\Gamma(Y)\rightarrow \Gamma(Z)$ of the
respective coordinate algebras. If, in addition, $Y \varsubsetneq
Z$ then $\ker \varphi \ne 1$.
\end{cor}

\section{Prevarieties} \label{sec:I-pvar}

We have already mentioned in Section \ref{sec:I-eofag} that
coordinate algebras of algebraic sets map into a cartesian power
of $B$. This observation suggests that this object deserves a
detailed study and, indeed, yields results about coordinate
algebras. In this section, we develop this approach. For that
purpose, we introduce the concept of  a prevariety and study
connections between prevarieties and certain classes of coordinate
algebras.

Given a class $\K$ of $A$-Lie algebras over a field $k$, we denote
by $S_A({\K})$ and $P_A({\K})$,  correspondingly, the classes of
all $A$-Lie subalgebras and all unrestricted $A$-cartesian
products of algebras from $\K$.

\begin{defn} \label{defn:pvar}
Let $\K$ be a class of $A$-Lie algebras over a field $k$. The
class $\K$ is called an \emph{$A$-prevariety} {\rm(}or, for
brevity, just prevariety{\rm)} if $\K = S_A P_A (\K)$.
\end{defn}

Assume that $\K$ is a class of $A$-Lie algebras over a field $k$.
Then the least  $A$-prevariety containing  $\K$ is called the
\emph{$A$-prevariety generated by the class $\K$} and is denoted
by {\rm $A  \pvar(\K)$}. It is easy to see that $A  \pvar(\K) =
S_A P_A (\K)$.

Observe that definition of the radical of an algebraic set can be
given in terms of intersections of kernels of  a certain
collection of homomorphisms; restricting this collection to
homomorphisms with images in a particular class of algebras, we
arrive at the concept  of a \emph{radical with respect to a
class}.

\begin{defn} \label{defn:I-radtoclass}
Let $\K$ be a class of $A$-Lie algebras over $k$, $C$ an  $A$-Lie
algebra and $S$ a subset of $C$. Consider the family of all
$A$-homomorphisms $\varphi _i: C \rightarrow D$ with $D \in \K$
and such that $\ker(\varphi _i) \supseteq S$. The intersection of
their kernels is called the \emph{radical of the set $S$ with
respect to the class $\K$},
$$ \rad_\K (S) = \rad_\K (S,C) =
\mathop \bigcap \limits_{S \subseteq \ker(\varphi _i )}
\ker(\varphi _i).$$
\end{defn}

Notice that if  $C=\AX$ and $\K=\left\{B\right\}$ then $\rad_B (S)
= \rad_\K (S,C)$.

\begin{lem} \label{lem:I-rad2}
Let $\K$ be a class of $A$-Lie algebras, $B$ an arbitrary $A$-Lie
algebra and $S \subseteq B$. Then
\begin{enumerate}
    \item {\rm$\rad_{\K} (S) \supseteq \verb"id" \left< S \right> $.}
    \item {\rm${\raise0.7ex\hbox{$B$} \!\mathord{\left/
 {\vphantom {B {\rad_{\K} (S)}}}\right.\kern-\nulldelimiterspace}
\!\lower0.7ex\hbox{${\rad_{\K} (S)}$}} \in A  \pvar(\K)$}.
    \item {\rm $\rad_{\K} (S,B)$}
 is the smallest ideal $I$ of the algebra $B$ containing $S$ and
such that \;{\rm ${\raise0.7ex\hbox{$B$} \!\mathord{\left/
 {\vphantom {B {I}}}\right.\kern-\nulldelimiterspace}
\!\lower0.7ex\hbox{${I}$}} \in A  \pvar(\K)$}.
    \item {\rm$\rad_{\K} (S) = \rad_{A  \pvar(\K)} (S)$}.
\end{enumerate}
\end{lem}

\begin{proof} (1) The first statement is straightforward.

(2) By Remak's theorem ${\raise0.7ex\hbox{$B$} \!\mathord{\left/
 {\vphantom {B {\rad_\K (S)}}}\right.\kern-\nulldelimiterspace}
\!\lower0.7ex\hbox{${\rad_\K} (S)$}}$
 $A$-embeds into the cartesian power $\prod\limits_{i \in I}
{{\raise0.7ex\hbox{$B$} \!\mathord{\left/
 {\vphantom {B {\ker \;\varphi _i }}}\right.\kern-\nulldelimiterspace}
\!\lower0.7ex\hbox{${\ker \;\varphi _i }$}}}$. For each $i \in I$
the algebra ${\raise0.7ex\hbox{$B$} \!\mathord{\left/ {\vphantom
{B{\ker \;\varphi _i }}}\right.\kern-\nulldelimiterspace}
\!\lower0.7ex\hbox{${\ker \;\varphi _i }$}}$ $A$-embeds into the
algebra $D_{\varphi _i }  \in \K$. Therefore the algebra
${\raise0.7ex\hbox{$B$} \!\mathord{\left/
 {\vphantom {B {\rad_{\K} (S)}}}\right.\kern-\nulldelimiterspace}
\!\lower0.7ex\hbox{${\rad_{\K} (S)}$}}$ is an $A$-subalgebra of
$\prod\limits_{i \in I} {D_{\varphi _i }}$.

(3) Let $J$ be an ideal of $B$, $J \supseteq S$ and
${\raise0.7ex\hbox{$B$} \!\mathord{\left/
 {\vphantom {B J}}\right.\kern-\nulldelimiterspace}
\!\lower0.7ex\hbox{$J$}} \in A  \pvar({\K})$. Then
${\raise0.7ex\hbox{$B$} \!\mathord{\left/
 {\vphantom {B J}}\right.\kern-\nulldelimiterspace}
\!\lower0.7ex\hbox{$J$}}$ is an $A$-subalgebra of $\prod\limits_{i
\in I} {D_i } $, where $D_i  \in {\K}$. Consequently, for every $i
\in I$ there exists an $A$-homomorphism $\varphi _i :B \to D_i $
such that $J \subseteq \ker \;\varphi _i $. Furthermore, $J =
\bigcap\limits_{i \in I} {\ker \;\varphi _i } $ and thus $J
\supseteq \rad_{\K} (S)$.

(4) The inclusion $\rad_{\K} (S) \supseteq \rad_{A  \pvar({\K})}
(S)$ is obvious. Suppose that $\rad_{\K} (S) \supsetneq \rad_{A
\pvar({\K})} (S)$ then  ${\raise0.7ex\hbox{$B$} \!\mathord{\left/
 {\vphantom {B {\rad_{A \pvar(\K)} (S)}}}\right.\kern-\nulldelimiterspace}
\!\lower0.7ex\hbox{$\rad_{A\pvar(\K) (S)}$}}  \in A \pvar({\K})$,
a contradiction with statement (3) of the lemma.
\end{proof}

Since the concept of a prevariety is extremely important, it will
be useful to have an alternative definition.

An $A$-Lie algebra $B$ is called \emph{$A$-approximated} by a
class $\K$, if for any $b \in B$, $b \ne 0$ there exists an
$A$-homomorphism $\varphi _b: B \rightarrow C$ for some $C \in \K$
and such that $\varphi _b (b) \ne 0$. The set of all $A$-Lie
algebras that are $A$-approximated by the class $\K$ is denoted by
\/ $\verb"Res"_A (\K)$.

\begin{lem} \label{lem:I-pvar=res}
For any class of $A$-Lie algebras $\K$  {\rm
\[
A  \pvar(\K) = \verb"Res"_A (\K).
\]
}
\end{lem}

\begin{proof} Clearly $P_A (\K) \subseteq \verb"Res"_A (\K)$, therefore $A
\pvar(\K) \subseteq \verb"Res"_A (\K)$. To prove the converse take
an arbitrary $A$-Lie algebra $B \in \verb"Res"_A (\K)$. According
to Definition \ref{defn:I-radtoclass}, $\rad_{\K}
(\left\{0\right\}, B) = 0$. Therefore $B \in A  \pvar(\K)$.
\end{proof}

Another very important feature of prevarieties is that every
prevariety has a theory of generators and relations.

\begin{lem} \label{lem:prpvar} Let $\K$ be
an $A$-prevariety of $A$-Lie algebras over a field $k$. Let $B$ be
an $A$-Lie algebra and  $B =  \left< X\mid R \right>$ its
presentation in the category of all $A$-Lie algebras, $R \subseteq
\AX$. Then the algebra $C$ lies in ${\K}$ if and only if\/
{\rm$\id \left< R \right>  = \rad_{\K} (R, \AX)$}.
\end{lem}

\begin{proof}
Let $\id \left< R \right>  =\rad_{\K} (R,\AX)$. Since $B \cong _A
{\raise0.7ex\hbox{$\AX$} \!\mathord{\left/
 {\vphantom {\AX \id \left< R \right>}}\right.\kern-\nulldelimiterspace}
\!\lower0.7ex\hbox{$\id \left< R \right>$}}$ then, according to
Lemma \ref{lem:I-rad2}, we have $B \in {\K}$. On the other hand
$\id \left< R \right> \subseteq \rad_{\K} (R,\AX)$. Therefore, if
$B \in {\K}$, then, by Lemma \ref{lem:I-rad2}, $\id \left< R
\right>  = \rad_{\K} (R,\AX)$.
\end{proof}

 By Lemma \ref{lem:prpvar}, for every algebra $B \in \K$
there exists a presentation $B = \left< X\mid R \right> _{\K}$, $R
\subseteq \AX$ of the form $$ B \cong _A
{\raise0.7ex\hbox{${\AX}$} \!\mathord{\left/
 {\vphantom {{\AX} {\rad_{\K} (R,\AX)}}}\right.\kern-\nulldelimiterspace}
\!\lower0.7ex\hbox{${\rad_{\K} (R,\AX)}$}}.$$ Indeed, if we treat
$B$ as an $A$-Lie algebra then $ B \cong _A
{\raise0.7ex\hbox{${\AX}$} \!\mathord{\left/
 {\vphantom {{\AX} {\texttt{id}\left< R \right> }}}\right.\kern-\nulldelimiterspace}
\!\lower0.7ex\hbox{${\texttt{id}\left< R \right>}$}}$. Next, since
$B \in \K$ we see that $\verb"id" \left< R \right> = \rad_{\K}
(R,\AX)$.

The algebra $$ A_{\K} \left[X\right] = {\raise0.7ex\hbox{${\AX}$}
\!\mathord{\left/
 {\vphantom {{\AX} {\rad_{\K} (0,\AX)}}}\right.\kern-\nulldelimiterspace}
\!\lower0.7ex\hbox{${\rad_{\K} (0,\AX)}$}}$$ is called \emph{the
free object of the prevariety $\K$}.

It might occur that $\id \left< R \right>$ is not finitely
generated, but $\rad_{\K} (R,\AX) =\rad_{\K} (R_0 ,\AX)$, where
$R_0 $ is a finite set. In which case the corresponding $A$-Lie
algebra $B$ is not finitely presented in the category of all
$A$-Lie algebras but is finitely presented in $\K$. Moreover, two
prevarieties $\K_1$ and $\K_2 $ coincide if and only we have, for
 any set $R \subseteq \AX$ and any set $X$,
 $$
\rad_{\K_1 } (R,\AX) = \rad_{\K_2 } (R,\AX).$$

We denote by $\K_\omega$ the class of all finitely generated
$A$-Lie subalgebras of $\K$.

\begin{lem} Let $S_A ({\K}) = {\K}$ and  $B$  a finitely
generated $A$-Lie algebra. Then {\rm$\rad_{\K} (R,B) =
\rad_{\K_\omega} (R,B)$} for any $R \subseteq B$.
\end{lem}

\begin{proof} Any $A$-homomorphism $\varphi:B \rightarrow D$, $D \in {\K}$
induces an $A$-homomorphism $\varphi _0 :B \rightarrow D_0$ such
that $D_0  = \im \;\varphi $ and $\ker \;\varphi  = \ker \;\varphi
_0 $. Since $D_0 $ is an $A$-subalgebra of  $D$, we have $K_0 \in
{\K}$. Finally, since $D_0  = \im \;\varphi $, the algebra $K_0 $
is a finitely generated $A$-algebra.
\end{proof}

The following lemma is a direct corollary of Lemma
\ref{lem:prpvar}.

\begin{lem} Let ${\K}_1 $ and ${\K}_2 $ be prevarieties of
$A$-Lie algebras. Then $({\K}_1 )_\omega   = ({\K}_2 )_\omega$ if
and only if for any finite set $X$ and any finite subset $R
\subseteq A[X]$ the radicals \, {\rm $\rad_{{\K}_1 } (R,\AX)$} and
\, {\rm $\rad_{{\K}_2 } (R,\AX)$} coincide.
\end{lem}

 The next theorem connects the theory of
prevarieties with algebraic geometry.

\begin{thm} \label{thm:I-pvar}
Let $B$ be an $A$-Lie algebra over a field $k$. Then all
coordinate algebras over $B$ lie in the prevariety {\rm $A
\pvar(B)$} and, conversely, every finitely generated $A$-Lie
algebra from the prevariety {\rm $A  \pvar(B)$} is a coordinate
algebra of an algebraic set over $B$.
\end{thm}

\begin{proof} Let $C=\left<X\mid S\right>_A$ be a finitely generated $A$-Lie algebra. By
Lemma \ref{lem:prpvar}, $C\in A  \pvar(B)$ if and only if $\id
\left< S \right>  = \rad_{A  \pvar(B)} (S, \AX)$. According to
Lemma \ref{lem:I-rad2}, $\rad_{A  \pvar(B)} (S) = \rad_{A
\pvar(B)} (S, \AX)$. Moreover,  $\rad_{B} (S) = \rad_{A \pvar(\K)}
(S, \AX)$. Consequently, $C\in A  \pvar(B)$  if and only if $\id
\left< S \right>= \rad_B(S)$. Therefore $C$ is a coordinate
algebra of an algebraic set of $V_B(S)$. \end{proof}

\begin{cor} \label{cor:I-coord=res}
A finitely generated $A$-Lie algebra $C$ is a coordinate algebra
over an $A$-Lie algebra $B$ if and only if $C$ is $A$-approximated
by $B$.
\end{cor}

\begin{rem}
In Section {\rm \ref{sec:eqthm}} we introduced the notion of a
co-presentation. We can  now reformulate this definition as
follows: a co-presentation $(X,S)$ of a coordinate algebra $\Gamma
_B (S)$ is its presentation in {\rm$A \pvar(B)$}.
\end{rem}

 If $C$ is a coordinate algebra of an irreducible
algebraic set over $B$ then it possesses yet another, stronger
property. It is very important in the study of universal closures
(see Sections \ref{sec:I-uncl} and \ref{sec:I-theunivcl}).

\begin{defn} \label{defn:I-discr}
An $A$-Lie algebra $C$ is said to be \emph{$A$-discriminated} by
an $A$-Lie algebra $B$ if for every finite subset $\left\{ c_1
,\ldots,c_m \right\}$ of nonzero elements from the algebra $C$
there exists an $A$-homomorphism $\varphi :C \rightarrow B$ such
that $\varphi (c_i ) \ne 0$, where $i = 1, \dots, m$. The set of
all $A$-Lie algebras discriminated by the algebra $B$ is denoted
by $\verb"Dis"_A (B)$.
\end{defn}

\begin{lem}
Let $B$ be an $A$-Lie algebra and  $C$ the coordinate algebra of
an irreducible algebraic set over $B$. Then $C \in \verb"Dis"_A
(B)$.
\end{lem}

\begin{proof} Let $C = \Gamma _B (Y) = \Gamma _B (S)$, where $Y$ is
an irreducible algebraic set over $B$. Assume that the statement
of the lemma is not true. Then there exists a tuple
\[
f_1  + \rad_B (S),\ldots,f_m  + \rad_B (S) \in \Gamma _B (S),\quad
f_i  \notin \rad_B (S),\;i = 1,\ldots,m,
\]
such that the image of at least one of these elements under every
$A$-homomorphism $\varphi: C \rightarrow B$ is zero. Let $Y_i$ be
the algebraic set over $B$ defined by  the system of equations $S
\cup \left\{ f_i \right\}$, $Y_i  = V_B (S \cup \left\{ f_i
\right\} )$. In this notation
\[
Y = Y_1  \cup \cdots \cup Y_m ,\quad Y \ne Y_i ,\quad i =
1,\dots,m.
\]
This leads to a contradiction, since $Y$ is an irreducible
algebraic set.
\end{proof}

 Let $\bar B = \prod\limits_{i \in I} {B^{(i)} } $ be
the unrestricted cartesian power of the algebra $B$. We turn the
algebra $\bar B$ into an $A$-Lie algebra by designating the
diagonal copy of the algebra $A$. Suppose that the cardinality of
the set $I$  is the maximum of the cardinalities of $B$ and
$\aleph_0$. Then  Theorem \ref{thm:I-pvar}, Proposition
\ref{prop:I-coormapstoB} and Definition \ref{defn:pvar} yield the
following result.

\begin{thm} \label{thm:I-Bn=coor}
Let $Y \subseteq B^n $ be an algebraic set over $B$. Then the
coordinate algebra $\Gamma (Y)$ $A$-embeds into the algebra $\bar
B$. Conversely, every finitely generated $A$-subalgebra of the
algebra $\bar B$ is the coordinate algebra of an algebraic set
over $B$.
\end{thm}

\begin{defn} \label{defn:realisation}
An $n$-generated $A$-subalgebra $C = \left< c_1 ,\ldots,c_n
\right>_A$ of  $\bar B$ is called the \emph{realisation of
$\Gamma_B (Y)$ in $\bar B$} if the complete set of relations $R
\subset \AX$ {\rm(}for the generators $c_1 ,\ldots,c_n  \in \bar
B${\rm)} coincides with {\rm$\rad(Y)$}.
\end{defn}

Note that in the above definition the number of generators $n$
coincides with the dimension of the affine space $B^n \supseteq Y$
and the generators are chosen in such a way that $R =\rad(Y)$.
However, the realisation of $\Gamma_B (Y)$ in  $\bar B$ is not
unique.

Let $C= \left< c_1 ,\ldots, c_n  \right>_A$ and $\tilde C= \left<
\tilde c_1 ,\ldots,\tilde c_n\right>_A$ be two $n$-generated
$A$-subalgebras of $\bar B$ and $R,\tilde R \subset \AX$ be the
respective complete sets of relations. The algebras $C$ and
$\tilde C$ are realisations of the coordinate algebra of the same
algebraic set if and only if $R = \tilde R$.

 Therefore Theorem \ref{thm:I-Bn=coor}  can  be refined as follows.

\begin{thm} \label{rem:I-diminB}
Let $Y \subseteq B^n $ be an algebraic set. Its coordinate algebra
$\Gamma_B (Y)$ has a realisation in $\bar B$, $\Gamma_B (Y) =
\left< c_1 ,\dots,c_n \right>_A$, $c_i \in \bar B$. Furthermore,
generators $c_1 ,\dots,c_n $ can be chosen in such a way that
$$Y = \left\{ (c_1 ^{(i)} ,\dots,c_n ^{(i)} )\;|\; i \in
I\right\}.$$

Conversely, if $C$ is a finitely generated $A$-subalgebra of $\bar
B$, $C =  \left< c_1 ,\dots,c_n  \right>_A$ then there exists a
unique algebraic set $Y$ such that $C$ is a realisation of\,
$\Gamma_B (Y)$ in $\bar B$. In this case $$Y \supseteq \left\{
(c_1 ^{(i)} ,\dots,c_n ^{(i)} )\; | \; i \in I\right\},;$$
moreover, the closure $$\overline{\left\{ (c_1 ^{(i)} ,\dots,c_n
^{(i)} )\; | \; i \in I\right\}}$$ in the Zariski topology
coincides with $Y$.
\end{thm}

 Theorem \ref{thm:I-pvar} demonstrates the importance of
the prevariety $A  \pvar(B)$ in studying algebraic geometry over
an $A$-Lie algebra. Unfortunately, the language of  prevarieties
is not very convenient, since prevarieties are not necessarily
axiomatisable classes. In the next section we introduce some
additional logical classes which have the advantage of being
axiomatisable, and discusst some of their connections with
prevarieties.

\section{Universal classes} \label{sec:I-uncl}

Given a class  $\K$  of  $A$-Lie algebras over a field $k$, we
construct several model-theoretical classes of $A$-Lie algebras.

First we need a number of definitions.

\bi \item An \emph{$A$-universal sentence} of the language $L_A$
is a formula of the type
\[
\forall x_1 \cdots\forall x_n \,\left(\mathop  \bigvee \limits_{j
= 1}^s \mathop  \bigwedge \limits_{i = 1}^t (u_{ij} (\bar x,\bar
a_{ij} ) = 0\; \wedge \;w_{ij} (\bar x,\bar c_{ij} ) \ne
0)\right),
\]
where $\bar x = (x_1 ,\ldots,x_n )$ is an $n$-tuple of variables,
$\bar a_{ij} $ and $\bar c_{ij} $ are sets of constants from the
algebra $A$ and $u_{ij}$, $w_{ij} $ are terms in the language
$L_A$ in variables $x_1 ,\ldots,x_n$. If an $A$-universal sentence
involves no constants from the algebra $A$ this notion specialises
into a standard notion of an universal sentence in the language
$L$.

 \item An \emph{$A$-identity} of the language $L_A$ is the formula of
the form
\[
\forall x_1 \cdots \forall x_n \left(\mathop  \bigwedge \limits_{i
= 1}^m \;r_i (\bar x, \bar a_{ij}) = 0\right),
\]
where $r_i (\bar x)$ are terms in the language $L_A$ in variables
$x_1 ,\dots,x_n$. If $A$-identity involves no constants from $A$
we come to the standard notion of identity of the language $L$.

 \item An \emph{$A$-quasi identity} of the language $L_A$  is a formula
of the form
\[
\forall x_1 \cdots\forall x_n \left(\mathop  \bigwedge \limits_{i
= 1}^m \;r_i (\bar x,\bar a_{ij}) = 0 \rightarrow s(\bar x, \bar
b) = 0\right),
\]
where $r_i (\bar x)$ and $s(\bar x)$ are terms. A
coefficients-free analogue is the notion of a quasi identity. \ei

Now we are ready to introduce the main definitions of this
section---they are standard concepts of universal algebra and
model theory.

 \bi \item A class of Lie algebras ${\K}$ is called a
\emph{variety} if it can be axiomatised by a set of identities.

\item A class of Lie algebras ${\K}$ is called a
\emph{quasivariety} if it can be axiomatised by a set of quasi
identities.

\item A class of Lie algebras ${\K}$ is called a \emph{universal
class} if it can be axiomatised by a set of universal sentences.
\ei

Replacing identity by $A$-identity, etc., we come to
``$A$-versions'' of these definitions: $A$-variety,
$A$-quasivariety and $A$-universal class.

\bi

\item The $A$-variety $A \var(\K)$ generated by the class $\K$ is
the class of all $A$-Lie algebras that satisfy all the identities
of the language $L_A$ satisfied by all algebras from $\K$.

\item The $A$-quasivariety $A \qvar(\K)$ generated by the class
$\K$ is the class of all $A$-Lie algebras that satisfy all the
quasi identities of the language $L_A$ satisfied by all algebras
from $\K$.

\item The $A$-universal closure $A  \fucl(\K)$ generated by the
class $\K$ is the class of all $A$-Lie algebras that satisfy all
the universal sentences of the language $L_A$ satisfied by all
algebras from  $\K$.

\ei

 For later use (see Section \ref{chap:II}) we also need to
consider the classes $\texttt{var}(\K)$, $\texttt{qvar}(\K)$ and
$\ucl(\K)$, which, by the definition, are the variety, the
quasivariety and the universal closure generated by the class $\K$
in the first order language $L$. The classes $\texttt{var}(\K)$,
$\texttt{qvar}(\K)$ and $\ucl(\K)$ are special cases of the
classes $A \var(\K)$, $A  \qvar(\K)$ and $A \fucl(\K)$.

Notice the following inclusions:
\[
A  \fucl(\K) \subseteq A  \qvar(\K) \subseteq A  \var(\K).
\]
The first inclusion is obvious, while the second one is implied by
the fact that every identity is equivalent to a conjunction of a
finite number of quasi identities. For instance, the identity
$$\forall x_1 \cdots \forall x_n \left(\mathop \bigwedge \limits_{i
= 1}^m \;r_i (\bar x) = 0\right)$$ is equivalent to the set of $m$
quasi identities $$\forall x_1 \cdots\forall x_n \forall y(y = y
\to r_i (\bar x) = 0).$$

Obviously the classes $A\var(\K)$ and  $A\qvar(\K)$ are
prevarieties. Moreover, $$A\pvar(\K) \subseteq A\qvar(\K)
\subseteq A\var(\K).$$

\section{Quasivarieties}

Every quasivariety is a prevariety and therefore quasivarieties
allow for a theory of generators and relations. In this section we
give a method of computation of radicals with respect to a
quasivariety.

Let $\K$ be a quasivariety of $A$-Lie algebras over a field $k$,
$B$ an $A$-Lie algebra over $k$ and $S$ a  subset of $B$.

Set $R_0$ to be the ideal of  $B$ generated by the set $S$ and
suppose that $R_i$ is already defined. Denote by $T_i$ the set of
all elements in  $B$ of the form $s(b_1,\dots,b_n )$, where $b_1,
\dots, b_n \in B$ and the quasi identity
\[
\forall x_1 \cdots\forall x_n \left(\mathop  \bigwedge \limits_{j
= 1}^m \;r_j (\bar x) = 0 \to s(\bar x) = 0\right)
\]
is true in all  algebras from $\K$, $r_j (b_1 ,\dots, b_n ) \in
R_i, \; j = 1, \dots, m$. In this notation we define the ideal
$R_{i + 1}=\id \left< R_i \cup T_i \right>$. As the result, we
have an ascending chain of ideals in $B$:
\[
R_0  \le \dots \le R_i  \le R_{i + 1}  \le \dots
\]

\begin{lem} \label{lem:I-radmeth}
In this notation, {\rm $\rad_{\K} (S,B) = \bigcup \limits_{i =
0}^\infty {R_i}$}.
\end{lem}

\begin{proof} Let $R = \bigcup\limits_{i = 0}^\infty {R_i }$.
An easy induction on $i$ shows that $\rad_{\K} (S,B) \supseteq R$.
To prove the reverse inclusion, it suffices to show $B/R \in \K$,
that is, verify that all  quasi identities which are true in every
algebra from the class $\K$ are true in the algebra $B/R$. The
latter easily follows from the definition of the set $R$.
\end{proof}

\begin{lem}
Let {\rm $f \in \rad_{\K} (S,B)$}. Then there exists a finite
subset $S_{0,f} \subseteq S$ such that {\rm $f \in \rad_\K
(S_{0,f} ,B)$}.
\end{lem}

\begin{proof} We use induction on $i$. If $f \in R_0$ then the statement
is obvious. Assume that it holds for all the elements from $R_i$
and consider an element $f \in R_{i + 1}$. By  definition of
$R_{i+1}$,  $f$ has the form $s(b_1, \dots,b_n )$. Using an
argument similar to that of Lemma \ref{lem:I-radmeth}, we have
$$s(b_1, \dots,b_n ) \in \rad_\K(\left\{r_j (b_1 ,\dots,b_n
),\;j=1,\dots m\right\}).$$ Since, by the inductive assumption,
there exist finite sets $S_j \subseteq S$ such that $r_j (b_1
,\dots,b_n ) \in \rad_{\K} (S_j )$, we have $S_{0,f} = S_1 \cup
\cdots \cup S_m $.
\end{proof}

As we mentioned earlier the prevariety $A\pvar(B)$ is of
exceptional importance for  algebraic geometry over the algebra
$B$. Nitice however that a prevariety is not in general an
axiomatisable class. On the other hand, every quasivariety admits
an axiomatic description and $A\pvar(B) \subseteq A\qvar(B)$. The
situation is clarified by  the following result of Malcev
\cite{Mal1}.

\begin{prop}
A prevariety is axiomatisable if and only if it is a quasivariety.
\end{prop}

This fact  leads  to the following question:
\begin{quote}
{\em For which $\K$ the prevariety\/ {\rm $\texttt{pvar}({\K})$}
is a quasivariety?}
\end{quote}

This question is known in the theory of quasivarieties as Malcev's
Problem \cite{Mal2}. Malcev himself gave the following sufficient
condition for the prevariety $\texttt{pvar}({\K})$ to be a
quasivariety \cite{Mal1}.

\begin{prop}
 Let ${\K}$ be an axiomatisable class of Lie algebras.
Then {\rm $\texttt{pvar}({\K})$} is a quasivariety.
\end{prop}

Recently Gorbunov \cite{Gor} found a complete solution of Malcev's
Problem.

Fortunately, if a system of equations $S(X) = 0$ has only finitely
many variables (as it is always the case in the algebraic geometry
over $B$) the radical $\rad_B(S)$ in $\AX$ depends only on
finitely generated $A$-Lie algebras from $A \pvar(B)$. So we can
refine the question.

\begin{quote}

 {\em The Restricted Malcev's Problem:} For which ${\K}$ the
subclasses of finitely generated $A$-Lie algebras in
$A\pvar({\K})$ and $A \qvar({\K})$ coincide? We are particularly
interested in $A\pvar({B})$ and $A \qvar({B})$.

\end{quote}

A criterion in  Lemma \ref{lem:I-qv=pv} gives an answer to this
question. For its formulation we need to define
$A$-$q_\omega$-compactness.

\begin{defn} \label{defn:I-qw}
An $A$-Lie algebra $B$ is called  \emph{$A$-\/$q_\omega$-compact}
if for any $n\in \mathbb{N}$, any system of equations $S \subseteq
A\left[x_1 ,\dots,x_n \right]$ and every its consequence {\rm $f
\in \rad_B (S)$} there exists a finite subsystem $S_{0,f}
\subseteq S$, $S_{0,f}  = \left\{ f_1 ,\dots,f_m \right\}$ such
that the following $A$-quasi identity is true in the algebra $B$:
\[
\forall x_1 \cdots\forall x_n \left(\mathop  \bigwedge \limits_{i
= 1}^m \;f_i (\bar x) = 0 \rightarrow f(\bar x) = 0\right),
\]
or, in other words,  {\rm $f \in \rad_B (S_{0,f})$}.
\end{defn}

The notion of a $q_\omega $-compact algebra is a natural
generalisation of  that of a Noetherian algebra. To say that an
algebra is Noetherian is equivalent to requiring that the finite
subsystem $S_0$ from  Definition \ref{defn:I-qw} is universal for
all the consequences of $S$. Also, $q_\omega$-compactness of an
algebra implies that this  finite subsystem $S_0$ is dependant on
the choice of a consequence of $S$.

\begin{rem}\label{rem:I-qw-noeth}
If an $A$-Lie algebra $B$ is  $A$-equationally Noetherian then $B$
is $A$-$q_\omega$-compact.
\end{rem}

Recall that a set of formulas $T$ is called \emph{compact} (for
the algebra $B$) if $B$ satisfies all formulas from $T$ whenever
every finite submodel  satisfies all formulas from $T$. The
expression `$q_\omega$-compactness' comes from the following
observation: an algebra $B$ is $A$-$q_\omega $-compact whenever
$B$ is compact with respect to the sentences of the language $L_A$
of the form
\[
T = \left\{ s = 0 \mid s \in S \right\}  \cup \left\{ f \ne
0\right\}.
\]

 B. Plotkin and others use different terminology for
saying that an algebra is $q_\omega$-compact or equationally
Noetherian: the corresponding expressions are \emph{logically
Noetherian} for $q_\omega$-compact algebraic structures and
\emph{algebraically Noetherian} for equationally Noetherian
structures.

\begin{expl} \label{ex:nonqwcomp}
Let $B$ be a nilpotent Lie algebra of class two given by the
following presentation in the variety of class $\leq 2$ nilpotent
Lie algebras:
\[
B = \langle a_i,b_i \ (i \in \mathbb{N}) \mid a_i \circ a_j = 0,
b_i\circ b_j = 0, a_i\circ b_j = 0 \;(i \neq j) \rangle.
\]
Then the infinite quasi identity
\[
\forall x \forall y \left(\bigwedge_{i \in \mathbb{N}}(x\circ a_i
= 0 \bigwedge_{j \in \mathbb{N}} x\circ b_j = 0) \rightarrow
x\circ y = 0\right)
\]
holds in $B$, but for any finite subsets $I,J$ of $\mathbb{N}$ the
following quasi identity does not hold in $B$:
\[
\forall x \forall y \left(\bigwedge_{i \in I}(x\circ a_i = 0
\bigwedge_{j \in J} x\circ b_j = 0) \rightarrow x\circ y =
0\right).
\]
Therefore the algebra $B$ is not $q_\omega$-compact.
\end{expl}

\begin{lem} \label{lem:I-qv=pv}
An $A$-Lie algebra $B$ is  $A$-$q_\omega$-compact if and only if
{\rm
\[
 (A  \pvar(B))_\omega   = (A  \qvar(B))_\omega  .
\]}
\end{lem}

\begin{proof} We observe first that $$(A  \pvar(B))_\omega   = (A
\qvar(B))_\omega$$ if and only if $$\rad_{A \pvar(B)} (S) =
\rad_{A \qvar(B)} (S)$$  for any system of equations $S \subseteq
\AX$ with a finite number of variables. Recall that, in view of
Lemma \ref{lem:I-rad2}, $\rad_{A  \pvar(B)} (S) = \rad_B (S)$.
Hence it suffices to show that the algebra $B$ is
$q_\omega$-compact if and only if  $\rad_{A  \qvar(B)} (S) =
\rad_B (S)$ for any system  $S \subseteq \AX$ with a finite number
of variables. Using Lemma \ref{lem:I-radmeth} and definition of
the sets $R_i, \; i\in \mathbb{N}$ we conclude that the
 $A$-$q_\omega$-compactness of $B$ follows from
the condition $\rad_{B} (S)=\rad_{A  \qvar(B)} (S)$ and hence from
the equality $\rad_{A \pvar(B)} (S) = \rad_{A \qvar(B)} (S)$.

Now we have to prove the reverse implication. Let $B$ be an
$A$-$q_\omega$-compact $A$-Lie algebra. We shall show that
$\rad_{A \qvar(B)} (S) = \rad_B (S)$ for an arbitrary system $S$
of equations. It follows from the definition of the radical with
respect to a quasivariety that $\rad_B (S) \supseteq \rad_{A
\qvar(B)} (S)$. Therefore it will suffice to prove that $f \in
\rad_{A \qvar(B)} (S)$ for any $f \in \rad_B (S)$. Since  $B$ is
$q_\omega $-compact then there exists a finite subsystem $S_{0,f}
\subseteq S$, $S_{0,f} = \left\{ f_1 ,\dots,f_m \right\}$ such
that the quasi identity
\[\forall x_1 \cdots\forall x_n \left(\mathop \bigwedge \limits_{j = 1}^m
\;f_j (\bar x) = 0 \rightarrow f(\bar x) = 0\right)\] holds in
$B$. Therefore this quasi identity is satisfied by every algebra
$C_i \in A  \qvar(B)$ and thus $f \in \rad_{A \qvar(B)} (S)$.
\end{proof}

\begin{cor} Let $B$ be an $A$-equationally Noetherian $A$-Lie algebra.
Then the class of all coordinate algebras of algebraic sets over
$B$ coincides with the class {\rm $(A  \qvar(B))_\omega$}.
\end{cor}

Myasnikov and Remeslennikov \cite[Problem 2 in Section 9]{alg2}
asked if there exists an $q_\omega$-compact groups which is not
equationally Noetherian; a similar question can be formulated for
an arbitrary algebraic structure.  The question was answered by
Goebel and Shelah \cite{Shelah} who constructed a group which is
not equationally Noetherian, but is $q_\omega$-compact; the same
construction works for Lie algebras. However, the author's
conversations  with B.~Plotkin, A.~Myasnikov and V.~Remeslennikov
led to an observation that all counterexamples known so far are
not finitely generated.

\begin{problem}
Are there  finitely generated $q_\omega$-compact $A$-Lie algebras
which are not $A$-equationally Noetherian?
\end{problem}

\section{Universal closure} \label{sec:I-theunivcl}

If an $A$-Lie algebra $B$ is $A$-equationally Noetherian then
every algebraic set over $B$ is a finite union of irreducible
algebraic sets and, moreover, this presentation is unique (Theorem
\ref{thm:1-1}). This shifts the focus of the study of algebraic
sets onto their irreducible components. It turns out that
irreducible algebraic sets and the corresponding coordinate
algebras are the algebraic counterparts to the universal closure
of the algebra $B$ (see Theorem \ref{thm:I-coor=ucl} below).

\begin{lem} \label{lem:I-dis-ucl}
Let $B$ and $C$ be $A$-Lie algebras such that
 $C \in \verb"Dis"_A (B)$ {\rm(}see Definition {\rm \ref{defn:I-discr})}.
Then  {\rm $C \in A  \fucl(B)$}.
\end{lem}

\begin{proof} Recall that the condition $C \in A  \fucl(B)$ means that every finite submodel of the algebra $C$
 $A$-embeds into the algebra $B$. This is obvious, since
 $C$ is $A$-discriminated by the algebra $B$.
\end{proof}

\begin{cor} \label{cor:I--inucl}
Let $B$ be an $A$-Lie algebra, and suppose that a finitely
generated $A$-Lie algebra $C$ is the coordinate algebra of an
irreducible algebraic set over $B$. Then {\rm $C \in A \fucl(B)$}.
\end{cor}

\begin{lem} \label{lem:I-dis}
Let $B$ be an $A$-equationally Noetherian $A$-Lie algebra and $C$
a finitely generated $A$-Lie algebra. If  {\rm $C \in A \fucl(B)$}
then $C \in \verb"Dis"_A (B)$.
\end{lem}

\begin{proof} Since $B$ is $A$-equationally Noetherian, we have, in view  of
Remark \ref{rem:I-qw-noeth} and Lemma \ref{lem:I-qv=pv}, that $$(A
\pvar(B))_\omega = (A  \qvar(B))_\omega.$$
 And since $A  \fucl(B) \subseteq A  \qvar(B)$
the algebra $C$ lies in the class $(A  \pvar(B))_\omega$, i.e. $C$
is a coordinate algebra of an algebraic set over $B$ and has the
form $C = \Gamma _B (S)$ for $S \subseteq \AX$, $X = \left\{ x_1,
\dots, x_n \right\}$. Since  $B$ is $A$-equationally Noetherian we
can assume without loss of generality that $S = \left\{ s_1
,\dots,s_p \right\}$ is a finite system of equations.

Assume now that the algebra $C$ is not discriminated by  $B$. Then
there exists an $m$-tuple
\[
f_1  + \rad_B (S),\dots,f_m  + \rad_B (S) \in \Gamma _B (S),\quad
f_i  \notin \rad_B (S),\;i = 1,\dots,m
\]
such that the image of at least one of these elements under every
$A$-homomorphism $\varphi: C \rightarrow B$ is zero. In other
words, this implies that the algebra $B$ satisfies the  universal
formula
\[
\Phi  = \forall x_1 \cdots\forall x_n \left(\mathop  \bigwedge
\limits_{j = 1}^p \;s_j (\bar x) = 0 \rightarrow \mathop  \bigvee
\limits_{i = 1}^m f_i (\bar x) = 0\right).
\]
This yields a contradiction, since $C \in A  \fucl(B)$ and the
formula $\Phi$ is not true in the algebra $C$. Indeed, after
substituting the elements \[x_1  + \rad_B (S), \dots, x_n  +
\rad_B (S) \in \Gamma _B (S)\] into $\Phi$ we see that the
consequence is a false statement while the condition is true.
\end{proof}

We denote by $\verb"LDis"_A (B)$ the set of all $A$-Lie algebras
such that every finitely generated $A$-Lie subalgebra of an
algebra from $\verb"LDis"_A (B)$ is $A$-discriminated by $B$. The
following result follows directly from the previous discussion.

\begin{thm} \label{thm:I-Ldis}
If $B$ is an $A$-equationally Noetherian $A$-Lie algebra then {\rm
\[
A  \fucl(B) = \verb"LDis"_A (B).
\]
}
\end{thm}

\begin{thm} \label{thm:I-coor=ucl}
Let $B$ be an $A$-equationally Noetherian $A$-Lie algebra and $C$
a finitely generated $A$-Lie algebra. Then $C$ is the coordinate
algebra of an irreducible algebraic set over $B$ if and only if\/
{\rm $C \in A
 \fucl(B)$}.
\end{thm}

\begin{proof} If $C$ is a coordinate algebra of an irreducible algebraic
set over $B$ then, in view of Corollary \ref{cor:I--inucl}, $C \in
A \fucl(B)$.

Now take a finitely generated $A$-Lie algebra $C$  from the class
$A \fucl(B)$. According to  Remark \ref{rem:I-qw-noeth} and Lemma
\ref{lem:I-qv=pv}, $(A  \pvar(B))_\omega = (A  \qvar(B))_\omega$.
Since $A  \fucl(B) \subseteq A  \qvar(B)$ the algebra $C$ lies in
the class $(A  \pvar(B))_\omega$, that is, $C$ is the coordinate
algebra of an algebraic set over $B$, $C = \Gamma _B (Y)$. By
virtue of Theorem \ref{thm:I-Ldis}, $C$ is $A$-discriminated by
the algebra $B$. We are left to show that the algebraic set $Y$ is
irreducible.

We argue towards a contradiction and assume that the algebraic set
$Y$ is reducible: $ Y = Y_1 \cup \cdots \cup Y_m$, $Y \ne Y_i$, $i
= 1, \dots, m$. Then
$$\rad_B (Y) = \rad_B (Y_1) \cap \cdots \cap \rad_B (Y_m ),$$
$\rad_B (Y) < \rad_B (Y_i)$, $i = 1,\dots,m$.

Let $f_i  \in \rad_B (Y_i) \smallsetminus \rad_B (Y), \; i =
1,\ldots,m$ and let $p\in Y$. For the time being we treat the
elements of $C$ as  polynomial functions (see Proposition
\ref{prop:I-coralg-pol}). Since every $A$-homomorphism $\varphi: C
\rightarrow B$ can be regarded as a substitution of a point $p \in
Y$ into the elements of $C$ (see Equation (\ref{eq:I-observ})),
the polynomial $f_i$ vanishes at   $p \in Y_i$. Since every point
of $Y$ is contained in at least one of its irreducible components,
we come to a contradiction with the fact that $C = \Gamma _B (Y)$
is $A$-discriminated by the algebra $B$.
\end{proof}

Theorem \ref{thm:I-Bn=coor} shows that  description of coordinate
algebras is equivalent to  description of finitely generated
subalgebras of $\bar B$. It turns out that any coordinate algebra
of an irreducible algebraic set embeds into the ultrapower
$\prod\limits_{i \in I} {B^{(i)} }/{\D}$ with respect to an
ultrafilter ${\D}$ over the set $I$ (for definition see
\cite{ChisRem,CK}). We turn the algebra $\prod\limits_{i \in I}
{B^{(i)} } / {\D}$ into an $A$-Lie algebra by designating the
diagonal copy of $A$. The following theorem follows from a
classical theorem by Malcev's  \cite{Mal2}.

\begin{thm} \

\begin{enumerate}
    \item Let $Y$ be an irreducible algebraic set over $B$.
Then the coordinate algebra $\Gamma_B (Y)$ $A$-embeds into an
ultrapower $\prod\limits_{i \in I} {B^{(i)} } /{\D}$, $|I| = |B|$.
    \item If $B$ is $A$-equationally Noetherian, then any
    finitely generated $A$-subalgebra of any ultrapower $ {B^J }/{\D}$
 is a coordinate algebra of an irreducible algebraic set over $B$.
\end{enumerate}
\end{thm}

\section{Geometric equivalence} \label{sec:GE}

B. Plotkin \cite{Plot1}  introduced an important notion of
geometrically equivalent algebraic structures. Myasnikov and
Remeslennikov \cite{alg2} discuss this notion in the case of
groups and observe  that all their results can be transferred to
an arbitrary algebraic structure. In this section, we transfer
their results to Lie algebras.

At the intuitive level of understanding, two Lie algebras are
geometrically equivalent if they produce identical algebraic
geometries.

$A$-Lie algebras $B$ and $C$ are called \emph{geometrically
equivalent} if for every positive integer $n$ and every system $S
\subseteq \AX$, $X = \{ x  _1 ,\dots,x_n \} $ the radicals \/ {\rm
$\rad_B (S)$} and \/ {\rm $\rad_C (S)$} coincide.

Since radicals completely determine  algebraic sets and their
coordinate algebras, the study of algebraic geometry over $B$ is
equivalent to the study of algebraic geometry over $C$, provided
that $B$ and $C$ are geometrically equivalent. In the latter case,
the respective categories of algebraic sets are isomorphic.

\begin{lem}\label{lem:I-paw=pbw}
$A$-Lie algebras $B$ and $C$ are geometrically equivalent if and
only if {\rm $$ (A  \pvar(B))_\omega   = (A  \pvar(C))_\omega.$$}
\end{lem}

\begin{proof} Assume that  $B$ and $C$ are geometrically equivalent. Then the families of coordinate algebras over $B$
and $C$ coincide and, consequently, the classes  $(A
\pvar(B))_\omega$ and $(A \pvar(C))_\omega$ also coincide by
virtue Theorem \ref{thm:I-pvar}.

Now suppose that $(A  \pvar(B))_\omega   = (A \pvar(C))_\omega$.
Let $S \subseteq \AX$ be an arbitrary system and consider its
radical $\rad_B (S)$. Applying Lemma \ref{lem:I-rad2}, we have $$
\rad_B (S) = \rad_{A  \pvar(B)} (S) = \rad_{A  \pvar(C)} (S) =
\rad_C (S). $$ By  definition, $B$ and $C$ are geometrically
equivalent.
\end{proof}

\begin{cor} \label{cor:I-ge-qvar=}
If $A$-Lie algebras $B$ and $C$ are geometrically equivalent then
{\rm
\[
A  \qvar(B) = A  \qvar(C).
\]}
\end{cor}

\begin{proof} By definition,
\[
A  \qvar(B) = A  \qvar(A  \pvar(B))=A  \qvar((A \pvar(B))_\omega).
\]
Now, in view of Lemma \ref{lem:I-paw=pbw}, $A\qvar(B) = A
\qvar(C)$.
\end{proof}

Observe that  geometric equivalence of $A$-Lie algebras $B$ and
$C$ does not follow from the coincidence of the quasivarieties
generated by $B$ and $C$. As Lemma \ref{lem:I-qw-ge} demonstrates,
a counterexample can be found only for non $q_\omega$-compact
$A$-Lie algebras. Such counterexamples indeed exist, see Example
\ref{ex:nonqwcomp}. We refer to \cite{alg2} for more detail.

Nevertheless, for $q_\omega $-compact Lie algebras the coincidence
of quasivarieties is equivalent to geometric equivalence.

\begin{lem} \label{lem:I-qw-ge}
Let $B$ and $C$ be $q_\omega$-compact $A$-Lie algebras. Then $B$
and $C$ are geometrically equivalent if and only if {\rm $A
\qvar(B) = A  \qvar(C)$}.
\end{lem}
\begin{proof} By Corollary \ref{cor:I-ge-qvar=}, the geometric equivalence of
the algebras $B$ and $C$ yields  $A  \qvar(B) = A \qvar(C)$.

Since $B$ and $C$ are $q_\omega$-compact, then
\[
(A  \pvar(B))_\omega   = (A  \qvar(B))_\omega \hbox{ and }  (A
\pvar(C))_\omega   = (A  \qvar(C))_\omega
\]
by Lemma \ref{lem:I-qv=pv}. Therefore the families of coordinate
algebras over $B$ and $C$ coincide and the algebras $B$ and $C$
are geometrically equivalent. \end{proof}

\section{Algebraic geometry over free metabelian Lie algebra}
\label{chap:II}

The objective of this section is to give a brief account of recent
results in algebraic geometry over free metabelian Lie algebras.
We follow \cite{preprint,AGFMLA1,AGFMLA2}.

Throughout this section $\F$  will denote a free metabelian Lie
algebra over a field $k$. We also use notation $\F_r= \left< a_1
,\dots,a_r \right>$ for the free metabelian Lie algebra of rank
$r$ with the set of free generators (the free base) $a_1
,\dots,a_r$.

Recall that a metabelian Lie algebra is a Lie algebra $A$ which
satisfies the identity: $$(a \circ b) \circ (c \circ d) =0.$$ We
denote  by $\Fit(A)$  the Fitting radical of the algebra $A$. It
is well-known that $\Fit (A)$ has a natural structure of a module
over the ring of polynomials, see \cite{AGFMLA1}.

\subsection{The $\Delta$-localisation and the direct module
extension of the Fitting's radical}

In this subsection we introduce  $U$-algebras;  for any
$U$-algebra $A$ over $k$ we describe two extensions  of the
Fitting radical of $A$ (introduced in \cite{AGFMLA1}). These
constructions play an important role in the study of the universal
closure of the algebra $A$ (see Sections \ref{sec:I-uncl} and
\ref{sec:I-theunivcl}).

Following  {\rm \cite{AGFMLA1}}, we call a metabelian Lie algebra
$A$ over a field $k$ a \emph{ $U$-algebra} if
\begin{itemize}
    \item {\rm $\Fit(A)$} is abelian;
    \item {\rm $\Fit(A)$} is a torsion free module over the ring of
    polynomials.
\end{itemize}

Let  $A$ be a $U$-algebra over a field  $k$ and $\left\{ {z_\alpha
,\;\alpha  \in \Lambda } \right\}$  a maximal set of elements from
$A$ linearly independent modulo $\Fit (A)$. We denote by $V$ the
linear span of this set over $k$.

Let $\Delta  =  \left< x_\alpha ,\alpha \in \Lambda \right>$ be a
maximal ideal of the ring $R = k\left[x_\alpha ,\alpha  \in
\Lambda \right]$. Denote by  $R_\Delta$ the localisation of the
ring $R$ with respect to $\Delta $ and by $\Fit_\Delta (A)$ the
localisation of the module $\Fit(A)$ with the respect to the ideal
$\Delta $, that is, the closure of $\Fit(A)$ under the action of
the elements of $R_\Delta$ (for definitions see \cite{Bourb} and
\cite{Lang}). Consider the direct sum $V \oplus \Fit_\Delta (A)$
of vector spaces over $k$. By definition, the multiplication on
$V$ is inherited from $A$ and the multiplication on $\Fit_\Delta
(A)$ is trivial. Set
\[
u \circ z_\alpha = u \cdot x_\alpha  ,\quad u \in \Fit_\Delta
(A),\;z_\alpha \in V,\quad u \cdot x_\alpha   \in \Fit_\Delta (A).
\]
and extend multiplication on $\Fit_\Delta (A)$ to  $V$ by
linearity.

In this notation, the algebra $A_\Delta$ is called the
\emph{$\Delta$-localisation} of the algebra $A$.

\begin{prop} Let $A$ be a $U$-algebra and   $A_\Delta  $
 its $\Delta $-localisation. Then the algebras $A$ and
$A_\Delta$ are universally equivalent, {\rm $\ucl(A) =
\ucl(A_\Delta )$}.
\end{prop}

Let $M$ be a torsion free module over the ring of polynomials $R$.
By definition, algebra $A \oplus M$ is the direct sum of
$k$-vector spaces $V \oplus \Fit(A) \oplus M$. To define the
structure of an algebra on $V \oplus \Fit(A) \oplus M$ we need to
introduce multiplication. The multiplication on  $V$ is inherited
from $A$. The multiplication on $\Fit(A) \oplus M$ is trivial.  If
$b \in M$ and $z_\alpha \in V$, we set $b \circ z_\alpha = b \cdot
x_\alpha$ and extend this  multiplication by $b$ to elements from
$V$ using the distributivity law. This operation turns $A \oplus
M$ into a metabelian Lie algebra over $k$.

\begin{prop} \label{prop:II:ucldirext}
Let  $A$ be a $U$-algebra and  $M$  a finitely generated module
over $R$. Then {\rm
\begin{itemize}
    \item  $\ucl(A) = \ucl(A \oplus M)$;
    \item  $A\fucl(A) = A \fucl (A\oplus M)$;
    \item  $\Fit (A \oplus M)=\Fit(A) \oplus M$.
\end{itemize}}
\end{prop}
The operations of $\Delta$-localisation and direct module
extension commute. Following \cite{AGFMLA2}, we denote by
$\F_{r,s}$ the direct module extension of the Fitting radical of
$\F_r$ by the free $R$-module $T_s $ of the rank $s$,
$\F_{r,s}=\F_r\oplus T_s$.

\subsection{The case of a finite field} \label{sec:II-1}

Our aim is to construct diophantine algebraic geometry over the
free metabelian Lie algebra $\F_r$. When $r=1$, algebraic geometry
over $\F_r$ degenerates to Example \ref{expl:I-dim1}. Therefore we
restrict ourselves to the non-degenerated case $r \ge 2$.

\subsubsection{The universal axioms for the $\Phi_r$-algebras}
\label{ss:II-univax-PhiR}

In this section we formulate two collections of universal  axioms
$\Phi _r$ and $\Phi '_r$ in the languages $L$ and $L_{\F_r}$,
which axiomatise the universal classes $\ucl(\F_r)$ and
$\F_r\fucl(\F_r)$. Most of these formulas are in the first order
language $L$ and consequently belong to both $\Phi _r $ and $\Phi
'_r $. We shall  write the two  series of formulas simultaneously,
at every step indicating the differences between $\Phi _r $ and
$\Phi '_r $.

Axiom $\Phi 1$ below is the metabelian identity, axiom $\Phi 2$
postulates that there are no non-abelian nilpotent subalgebras and
axiom $\Phi 3$ is the $CT$-axiom:
\begin{gather}\notag
\begin{split}
\Phi 1: & \ \forall x_1 ,x_2 ,x_3 ,x_4 \quad (x_1 x_2 )(x_3 x_4 )
=0. \\
\Phi 2: & \ \forall x\forall y\quad xyx = 0 \wedge xyy = 0 \to xy
=0.\\
\Phi 3: & \ \forall x\forall y\forall z\quad x \ne 0 \wedge xy = 0
\wedge xz = 0 \to yz = 0.
\end{split}
\end{gather}
We next introduce an universal formula $\sFit(x)$ of the language
$L$ and in one variable
 which defines the Fitting radical
\begin{equation} \label{eq:fit1}
\sFit(x) \equiv (\forall y\;\,xyx = 0).
\end{equation}
An analogue of formula (\ref{eq:fit1}) in the language $L_{\F_r }$
is
\begin{equation} \label{eq:fit2}
\sFit '(x) \equiv \,\mathop  \bigwedge \limits_i \,(xa_i x = 0)\,
.
\end{equation}

From now on we restrict ourselves to the case of a finite field
$k$. In particular, the vector space
 {\rm ${\raise0.7ex\hbox{${\F_r }$} \!\mathord{\left/
 {\vphantom {{\F_r } {\Fit(\F_r )}}}\right.\kern-\nulldelimiterspace}
\!\lower0.7ex\hbox{${\Fit(\F_r )}$}}$} is finite and its dimension
over $k$ is $r$.

\begin{lem} \label{lem:II-linindmodfit}
Let $k$ be a finite field and $n \in \mathbb{N}$, $n \le r$. Then
the formula {\rm
\[
\varphi (x_1 ,\ldots,x_n) \equiv \mathop \bigwedge \limits_{(\alpha
_1 ,\ldots,\alpha _n ) \ne \bar 0} \neg \sFit (\alpha _1 x_1  +
\cdots + \alpha _n x_n )
\]}
of the language
$L$ is true on the elements $\left\{ b_1 ,\ldots,b_n \right\}$ of
$\F_r$ if and only if \/ $b_1 ,\ldots,b_n$ are linearly
independent modulo {\rm $\Fit(\F_r )$}.
\end{lem}

We can express the dimension axiom using formula $ \varphi$:
\[
\Phi 4: \forall x_1 \cdots \forall x_{r + 1} \;\neg \varphi (x_1
,\ldots,x_{r + 1} ).
\]
Since $\varphi $ is an existential formula, the formula $\neg
\varphi $ and  hence axiom $\Phi 4$ are universal. This axiom
postulates that the dimension of the factor-space
${\raise0.7ex\hbox{${B}$} \!\mathord{\left/
 {\vphantom {{B } {\Fit(B )}}}\right.\kern-\nulldelimiterspace}
\!\lower0.7ex\hbox{${\Fit(B)}$}}$ is at most $r$, provided that
$B$ satisfies $\Phi 1- \Phi 3$.

Recall that  $\Fit(\F_r )$ has a structure of a module over the
ring of polynomials $R = k\left[x_1 ,\dots,x_r \right]$. The
series of axioms $\Phi 5$, $\Phi' 5$,  $\Phi 6$, $\Phi 7$ and
$\Phi' 7$ express module properties of  $\Fit(\F_r )$. In this
axioms we use module (over $R$) notation. In particular, this
involves  rewriting  the polynomial $f(x_1 ,\ldots,x_n)$, $n \le
r$  in the signature of metabelian Lie algebra see \cite{AGFMLA1}
for a more detailed description of transcription from the module
signature to the Lie algebra signature.

The Fitting radical of the free metabelian Lie algebra is a
torsion free module  over the ring $R$. We use this observation to
write the following infinite series of axioms. For every nonzero
polynomial $f \in k \left[x_1 ,\ldots,x_n \right]$, $n \le r$,
write
\begin{gather} \notag
\begin{split}
\Phi 5:\forall z_1 \forall z_2\forall x_1 \cdots\forall x_n (z_1
z_2
\cdot f(x_1 ,\dots,x_n )  = 0 \, \wedge \, z_1 z_2  \ne 0)\\
\rightarrow (\neg \, \varphi (x_1 ,\ldots,x_n )).
\end{split}
\end{gather}
\normalsize Since  $\varphi (x_1 ,\dots ,x_n )$ is a
$\exists$-formula the formula $\Phi 5$ is a $\forall$-formula.

This property can be expressed in the language $L_{\F_r}$ as
follows. For any nonzero polynomial $f \in k \left[x_1 ,\ldots,x_r
\right]$ write:
\[
\Phi' 5: \ \forall z_1 \forall z_2 \quad (z_1 z_2 \cdot f(a_1
,\ldots,a_r ) = 0 \to z_1 z_2  = 0).
\]
The main advantage of this formula is that it does not involve the
formula $\varphi$, thus the restriction on the cardinality of the
field $k$ is not significant.

For every nonzero Lie polynomial $l(a_1 ,\ldots,a_n )$, $n \le r$
in variables $a_1 ,\ldots,a_r $ from the free base of $\F_r$, we
write
\[
\Phi 6: \ \forall x_1 \cdots \forall x_n \quad \varphi (x_1
,\ldots,x_n ) \to (l(x_1 ,\dots,x_n ) \ne 0).
\]
Since  $\varphi (x_1 ,\ldots ,x_n )$ is a $\exists$-formula the
formula $\Phi 6$ is a $\forall$-formula.

Series of axioms $\Phi 7$ and $\Phi'7$ are quite sophisticated. We
first introduce higher-dimensional analogues of Formulas
(\ref{eq:fit1}) and (\ref{eq:fit2}):
\begin{gather} \notag
\begin{split}
& \sFit(y_1,\ldots,y_l ;\,x_1 ,\ldots,x_n ) \equiv \left(\,\mathop
\bigwedge \limits_{i = 1}^n \,(y_1 x_i y_1 = 0)\,) \wedge \cdots
\wedge (\,\mathop  \bigwedge \limits_{i = 1}^n \,(y_l x_i y_l  =
0)\right), \\
& \sFit'(y_1 ,\ldots,y_l ) \equiv \sFit'(y_1 ) \wedge \cdots
\wedge \sFit'(y_l ).
\end{split}
\end{gather}

We begin with the series of axioms $\Phi' 7$ in the language
$L_{\F_r }$. Let $S$ be a fixed finite system of module equations
with variables $y_1 ,\ldots,y_l $ over the module $\Fit(\F_r )$.
Every equation from $S$ has the form
\[
h = y_1 f_1 (\bar x) + \cdots + y_l f_l (\bar x) - c = 0,\;\quad c
= c(a_1 ,\ldots,a_r ) \in \Fit(\F_r ),
\]
where $\bar x = \{ x_1 ,\ldots,x_r \} $ is a vector of variables
and $f_1 ,\ldots,f_l  \in R = k\left[\bar x\right]$. Suppose that
$S$ is inconsistent over  $\Fit(\F_r )$. This fact can be written,
in an obvious way, by logical formula in the signature of a
module. The system $S$ gives rise to a system of equations $S_1$
over $\F_r$. Replace every module equation $h_i = 0$ from $S$ by
the equation $h'_i  = 0$, $i = 1,\ldots,m$ in the signature of
$L_{\F_r}$ (see \cite{AGFMLA1}). This procedure results in a
system of equations $S_1$ over $\F_r$. For every inconsistent
system of module equations $S$, write
\[
\Phi '7: \ \psi '_S  \equiv \forall y_1 \cdots \forall y_l \quad
\sFit'(y_1 ,\dots,y_l ) \rightarrow \mathop  \bigvee \limits_{i =
1}^m \;h'_i (y_1 ,\ldots,y_l ) \ne 0.
\]
Notice that we have not used the restriction on the cardinality of
the field $k$.

Now we turn our attention  to the collection $\Phi_r$ in the
language $L$. Let $S$ be a system of $m$ module equations
inconsistent over the Fitting radical {\rm$\Fit_\Delta (\F_n )$}
of $\Delta $-local Lie algebra $(\F_n )_\Delta$. We write, for
every $n \in \mathbb{N}$, $n \le r$ and every such system $S$:
\begin{gather} \notag
\begin{split}
\Phi 7: \ \psi _{n,S}  \equiv \forall x_1 \cdots\forall x_n
\;\forall y_1 \cdots\forall y_l & \quad \varphi (x_1 ,\ldots,x_n )
\wedge
\sFit(y_1 ,\ldots,y_l ) \rightarrow \\
& \rightarrow \mathop \bigvee \limits_{i = 1}^m h_i (y_1
,\ldots,y_l ;\,x_1 ,\ldots,x_n ) \ne 0.
\end{split}
\end{gather}
The Lie polynomials $h_i$, $i = 1,\dots,m$ are constructed from
the system $S$ according to the following procedure. Consider the
$i$-th equation of $S$. It has the form
\begin{gather} \notag
\begin{split}
h'_i = y_1 f_1 (x_1 ,\ldots,x_n ) + \cdots + & y_l f_l (x_1
,\ldots,x_n) - c = 0,\\
& f_i  \in R,\;c = c(a_1 ,\ldots,a_n ) \in \Fit(\F_n ).
\end{split}
\end{gather}
After replace every occurrence of $a_j $ in $c(a_1 ,\ldots,a_n )$
by $x_j $, $j = 1,\ldots,n$ and rewriting the polynomial $h'_i$ in
the signature of Lie algebras (see \cite{AGFMLA1}), the resulting
Lie polynomial is $h_i$.

Denote, correspondingly, by $\Phi _r $ and by $\Phi '_r $ the
universal classes axiomatised by $\Phi 1$-$\Phi 7$ and by $\Phi
1$-$\Phi 4, \Phi' 5, \Phi 6, \Phi'7$. The algebras from these
classes are called, correspondingly, \emph{$\Phi _r $-algebras}
and \emph{$\Phi '_r$-algebras}.

\subsection{Main results}
Assume that the ground field $k$ is finite.

\begin{thm} \label{thm:II-ucl}
Let $A$ be an arbitrary finitely generated metabelian Lie algebra
over a finite field $k$. Then the following conditions are
equivalent:
\begin{itemize}
    \item {\rm$A \in \ucl(\F_r )$}.
    \item There exists $s \in \mathbb{N}$ so that
    $A$ is a subalgebra of $\F_{r,s}$.
    \item $A$ is a $\Phi_r$-algebra.
\end{itemize}
\end{thm}

\begin{cor} \label{cor:II-ucl}
The universal closure {\rm  $\ucl(\F_r )$} of the free metabelian
Lie algebra  $\F_r$ is axiomatised by $\Phi _r $.
\end{cor}

\begin{thm} \label{thm:II-F-ucl}
Let $A$ be an arbitrary finitely generated metabelian $\F_r$-Lie
algebra over a finite field $k$. Then the following conditions are
equivalent:
\begin{itemize}
    \item {\rm$A \in \F_r\fucl(\F_r )$};
    \item $A$ is a $\Phi'_r$-algebra;
    \item $A$ is $\F_r$-isomorphic to the algebra
     $\F_{r}\oplus M$, where $M$ is a torsion free module over
      $R=k\left[x_1 ,\dots,x_n \right]$.
\end{itemize}
\end{thm}

It is shown in {\rm\cite{AGFMLA1}}  that for every torsion free
$R$-module $M$ the algebra $\F_{r}\oplus M$ $\F_{r}$-embeds into
$\F_{r,s}$ for some $s \in \mathbb{N}$. It follows that all $\Phi'
_r$-algebras can be treated as $\F_{r}$-subalgebras of $\F_{r,s}$.

\begin{cor} \label{cor:II-F-ucl}
The universal closure  {\rm $\ucl(\F_r )$} of the free metabelian
Lie algebra $\F_r$ is axiomatised by $\Phi' _r $.
\end{cor}

\begin{thm}\label{thm:II-Phi} \label{thm:II-Phi'}
Axioms $\Phi _r $ {\rm (}and $\Phi' _r ${\rm)} form a recursive
set. The universal theory in the language $L$ {\rm (}in the
language $L_{\F_r}$, respectively{\rm)} of the algebra $\F_r$
{\rm(}treated as an $\F_r$-Lie algebra, respectively{\rm)} over a
finite field $k$ is decidable.
\end{thm}

\begin{thm}
Compatibility problem for a system of equations over the free
metabelian Lie algebra $\F_r$ is decidable.
\end{thm}

This result contrasts  with a result of Roman'kov \cite{Roman} on
the compatibility problem over  metabelian groups: he proved that
in the case of free metabelian groups of a sufficiently large rank
the problem is undecidable. His argument holds for free metabelian
Lie rings and free metabelian Lie algebras, provided that the
compatibility problem for the ground field is undecidable (and
thus the ground field is infinite).

We next classify all irreducible algebraic sets over $\F_r$.
Combining together Theorem  \ref{thm:I-coor=ucl} and
\ref{thm:II-F-ucl} we state

\begin{prop} \label{prop:II-irrcoralg} Let $\Gamma$ be an $\F_r$-Lie algebra.
Then $\Gamma$ is a coordinate algebra of an irreducible algebraic
set over $\F_r $ if and only if $\Gamma$ is $\F_r $-isomorphic to
$\F_r \oplus M$, where $M$ is a torsion free module over
$k\left[x_1 ,\ldots,x_r \right]$.
\end{prop}

Let $R = k\left[x_1,\ldots,x_r \right]$ and let $M$ be a finitely
generated torsion free module over $R$. Let $\hom_R (M,\Fit(\F_r
))$ be the set of all  $R$-homomorphisms from  $M$ to $\Fit(\F_r
)$ treated as a module over $R$. The lemma below describes the
canonical implementation of an algebraic set.

\begin{lem} In this notation, we have one-to-one correspondences
{\rm $$\hom_R (M,\Fit(\F_r )) \leftrightarrow \hom_{\F_r } (\F_r
\oplus M,\F_r ) \leftrightarrow Y,$$} where $Y$ is an irreducible
algebraic set over $\F_r $ such that $\Gamma(Y)=\F_r \oplus M$.
\end{lem}

\begin{thm}
Up to isomorphism, every irreducible algebraic set over $\F_r$ is
either {\rm $\hom_R(M,\Fit(\F_r ))$} for some finitely generated
torsion free module  $M$ over the ring $R$, or a point.
\end{thm}

\begin{cor}
Every irreducible algebraic set in the affine space $\F_r^n $, $n
= 1$, is, up to isomorphism, either a point or {\rm$\Fit(\F_r )$}.
\end{cor}

Recall that the rank $r(M)$ of the module $M$ over a ring $R$ is
the supremum of cardinalities of linearly independent over $R$
sets of elements from $M$. We set $r(\Gamma (Y))=r(M)$ if
$\Gamma(Y)=\F_r \oplus M$.

\begin{thm} \label{thm:II-dim}
For an  irreducible algebraic set $Y$ over $\F_r$
\[
\dim (Y) = r(\Gamma (Y)) = r(M).
\]
\end{thm}

\begin{rem}
When this paper had already been written E.~Daniyarova published her
PhD thesis, where, in particular, she obtains a description of the
quasivarieties {\rm$\texttt{qvar}(\F_r)$} and
{\rm$\F_r\qvar(\F_r)$}.
\end{rem}

\section{Algebraic geometry over a free Lie algebra}
\label{sec:AGOFLA}

In this section we outline results by Daniyarova and Remeslennikov
\cite{AGOFLA} on diophantine geometry over the free Lie algebra.

The aim of algebraic geometry is to classify irreducible algebraic
sets and their coordinate algebras. We expect that  the
classification problem of algebraic sets and coordinate algebras
in the case of free Lie algebra is very difficult if taken in its
full generality. We treat this problem only in two following
cases:
\begin{itemize}
    \item for algebraic sets defined by systems of equations with
    one variable;
    \item for bounded algebraic sets (see Definition
    \ref{defn:boundAS}).
\end{itemize}
In these cases, we reduce the problem to the corresponding problem
in the diophantine geometry over the ground field $k$.

The classification of algebraic sets and coordinate algebras over
free groups was considered in a series of papers
\cite{Appel,alg1,ChisRem, Guirardel,Lorents, Lorents1,
Lyndon,alg2}, the earliest of which, \cite{Lyndon}, dates back to
1960. In that paper Lyndon studied one-variable systems of
equations over free groups. Recently, a satisfactory
classification of irreducible algebraic sets over free group
$\FrGr$ and their coordinate groups  was given in  \cite{KM1,
KM2}. We begin this section by comparing results for free Lie
algebras \cite{AGOFLA} and  for free groups \cite{ChisRem}.

Let $\FrGr$ be a free group and $S$ any system of equations in one
variable over $\FrGr$ (i.e. $S\subseteq \FrGr \ast \langle
x\rangle$), such that $V_\FrGr(S)\ne\emptyset$. The full
description of all algebraic sets in $\FrGr^1=\FrGr$ is given by
the following two theorems.

\begin{thm} Any coordinate group $\Gamma_\FrGr(Y)$
of an irreducible algebraic set $Y\subseteq \FrGr^1$ satisfies one
of the following conditions.
\begin{itemize}
    \item $\Gamma_\FrGr(Y)\cong \FrGr$;
    \item $\Gamma_\FrGr(Y)\cong \FrGr\ast\langle x\rangle$, here $\ast$ is the sign of free product of groups;
    \item $\Gamma_\FrGr(Y)\cong \left< \FrGr,t \ |\; \left[u,t\right]=1\right>$, where $u$ generates a maximal
    cyclic subgroup of  \/ $\FrGr$.
\end{itemize}
\end{thm}

\begin{thm}
If\/ $V\ne\FrGr^1$ is an irreducible algebraic set defined by a
system of equations with one variable then:
\begin{itemize}
    \item either $V$ is a point or
    \item there exist elements $f$, $g$, $h\in \FrGr$ such that $V=fC_\FrGr(g)h$, where $C_\FrGr(g)$ stands
    for the centraliser of $g$ in $\FrGr$.
\end{itemize}
\end{thm}

On the other hand, the classification of algebraic sets in the
case of one-variable equations over free Lie algebra is much more
complicated.

\begin{thm} \label{thm:1var}
Up to isomorphism, an algebraic set defined by a consistent system
of equations with one variable over $F$ is either bounded
{\rm(}see Definition {\rm\ref{defn:boundAS})} or coincides with
$F$.

\end{thm}

Roughly speaking, bounded algebraic sets are the ones contained in
a finitely dimensional affine subspace of $F$. To make this
definition explicit, we introduce   the notion of a
parallelepipedon (Section \ref{ss:III-findimalgset}) and then show
that algebraic geometry in a parallelepipedon is equivalent to
diophantine geometry over the ground field (Section
\ref{ss:interp}).

\subsection{Parallelepipedons}
\label{ss:III-findimalgset}

We shall now show that every finite dimensional affine subspace
$V$ of a finitely generated free Lie algebra $F$ is an algebraic
set.

Let $V$ be a finite dimensional subspace of $F$, $V=\verb"lin"_k
(v_1 ,\ldots, v_m )$ where $\verb"lin"_k$ is the linear span over
$k$. Set
\[
s_1 (x) = x \circ v_1; \quad s_2 (x) = s_1 (x,v_1 ) \circ s_1 (v_2
,v_1 ) = (x \circ v_1 ) \circ (v_2 \circ v_1 ).
\]
Recursively, set
\begin{equation}  \label{eq:III-s_k}
 s_m (x)  =  s_{m - 1} (x,v_1 ,\ldots,v_{m - 1} ) \circ s_{m - 1} (v_m, v_1
,\ldots,v_{m - 1}).
\end{equation}
Notice that equations $s_l (x) = 0$ are linear over  $k$ and thus
define vector subspaces of the algebra $F$.

\begin{prop} \label{prop:III-lin-alg}
Any finite dimensional linear subspace of the algebra  $F$ is an
algebraic set over $F$. Furthermore, an $l$-dimensional subspace
of $F$ can be defined by an equation of the type $s_l$.
\end{prop}

 Moreover, affine shifts of linear subspaces (affine
subspaces) of $F$ are also algebraic sets. Let $c \in F$ be an
arbitrary element. Then the subspace $V + c$ (here $V$ is a linear
subspace of dimension $l$) is, obviously,  also an algebraic set
over the free Lie algebra $F$ in the obvious way: it is  defined
by the equation $s_l (x - c) = 0$.

\begin{cor}
Every affine shift of an arbitrary finite dimensional linear
subspace of the algebra $F$ is an algebraic set over $F$.
\end{cor}

 Using the same argument for systems of equations with
$n$ variables we get

\begin{cor} \label{cor:III-prodaff-alg}
Let  $V_i$, $i = 1,\dots,n$, be finite dimensional linear {\rm(}or
affine{\rm)} subspaces of the free Lie algebra $F$. The direct
product $\V = V_1 \times \cdots \times V_n \subset F^n$ is an
algebraic set over $F$.
\end{cor}

We denote  by $\V$ the direct product of affine finite dimensional
subspaces of $F$, $$\V =(V_1+c_1)\times \cdots \times (V_n+c_n)
\subset F^n,$$ $V_i=\verb"lin"_k\left\{v_1^i, \dots,
v_m^i\right\}.$ We call such spaces \emph{parallelepipedons}.

Let $\bar F=\prod\limits_{i \in I} {F^{(i)} }$ and $\bar k =
\prod\limits_{i \in I} {k^{(i)} }$ be the unrestricted cartesian
powers of $F$ and $k$, where the cardinality of the set $I$
coincides with the cardinality of $F$. We turn the algebra $\bar
F$ into an $F$-Lie algebra by identifying the diagonal copy of the
algebra $F$ with $F$. We next apply Theorem \ref{rem:I-diminB} and
obtain

\begin{prop} \label{prop:III-coordalgofaff}
Let $\Gamma (\V)=\left<\x_1 ,\ldots,\x_n  \right>_F$ be a
realisation {\rm(}see Definition {\rm\ref{defn:realisation})} of\/
$\Gamma (\V)$ in $\bar F$, $\x_1 ,\ldots,\x_n \in \bar F$. Then
the generators $\x_i$'s have the form
\begin{equation} \label{eq:xi}
\x_i  = t_1^i v_1^i + \cdots + t_{m_i }^i v_{m_i }^i ,\quad t_1^i
,\dots,t_{m_i }^i \in \bar k.
\end{equation}
and the coefficients $t_j^i$'s satisfy the conditions
\begin{itemize}
    \item $\left\{({t_j^i}^{(l)}),\ l\in I\right\}=k^M$ if the field $k$ is finite,
    \item the elements $t_i^j$'s are algebraically independent
    over $k$, otherwise.
\end{itemize}
\end{prop}

\subsection{Bounded algebraic sets and coordinate algebras}

\begin{defn} \label{defn:boundAS}
An algebraic set $Y$ over $F$ is called \emph{bounded} if it is
contained in a parallelepipedon.
\end{defn}

We list some  properties of bounded algebraic sets:

\begin{itemize}
    \item An arbitrary bounded algebraic set $Y$ is contained in
    infinitely many parallelepipedons.
    \item Let $Y, Z \subset F^n$ be algebraic sets and let $Y$ be
    bounded. Suppose next that there exists an epimorphism $\phi: Y
    \rightarrow Z$. Then $Z$ is also bounded.
\end{itemize}

This section is mostly concerned with the parallel concept  of
bounded coordinate algebra.

Set
\[
B(\bar F)= \left\{\xi \in \bar F| \hbox{ the degrees of the
coordinates of } \xi
 \hbox{ are bounded above} \right\}.
\]

\begin{prop} \label{prop:III--B(F)=Bok} In this
notation, $B(\bar F)\cong F \otimes_k \bar k$.
\end{prop}

A coordinate algebra is called
\emph{bounded} if it has a realisation in $B(\bar F)$ {\rm(}see
Definition {\rm\ref{defn:realisation})}.

\begin{prop} \label{prop:boundsA} \
\begin{itemize}
\item[1.] An algebraic set $Y \subseteq F^n$ over $F$ is bounded
if and only if its coordinate algebra $\Gamma_F(Y)$ is bounded.

\item[2.] If $Y \subseteq \V$ then any realisation of
$\Gamma_F(Y)=\left<\x_1, \ldots, \x_n\right>$ has the form
{\rm(\ref{eq:xi})}.

\item[3.] If the generators $\x_i$'s in a realisation of $\Gamma_F
(Y)$
 have the form {\rm(\ref{eq:xi})} then $Y \subset \V$.
\end{itemize}
\end{prop}

\subsection{The correspondences between algebraic sets, radicals and
coordinate algebras} \label{ss:interp}

The objective of this  Section is to establish a correspondence
between bounded algebraic geometry over the free Lie algebra and
algebraic geometry over the ground field $k$.

We begin with bounded algebraic geometry in dimension one, that
is, consider algebraic sets $Y \subset F$. Fix a parallelepipidon
$\V =V+c \subset F$, where $V=\verb"lin"_k\left\{v_1, \ldots,v_m
\right\}$, $\dim V=m$, $c\in F$. We shall show that there exists a
one-to-one correspondence between bounded algebraic sets over $F$
from $\V$ and algebraic sets over $k$ that lie in the affine space
$k^m$: $$Y_F \subseteq \V \leftrightarrow Y_k \subseteq k^m.$$ Our
correspondence takes the following form. Let $Y_F \subseteq \V$
then
\begin{equation} \label{eq:ktoF}
Y_k=\left\{( \alpha_1, \dots, \alpha_m) \in k^m\; | \; \alpha_1
v_1 + \cdots + \alpha_m v_m+c \in Y_F \right\}, Y_k \subseteq k^m.
\end{equation}
Conversely, let $Y_k$ be an algebraic set in $k^m$. Set
\begin{equation} \label{eq:Ftok}
Y_F=\left\{ \alpha_1 v_1 + \cdots +\alpha_m v_m\;|\; ( \alpha_1,
\dots, \alpha_m) \in Y_k\right\}, Y_F \subseteq \V.
\end{equation}
We wish to show that sets defined by (\ref{eq:ktoF}) and
(\ref{eq:Ftok}) are algebraic. We do this constructing the
corresponding systems of equations, $Y_F=V(S_F)$ and $Y_k=V(S_k)$.

As shown in Lemma \ref{lem:rad}, an algebraic set is uniquely
determined by its radical, and it will more convenient for us to
work with. First we define the correspondence between polynomials
in $\Fx$ and $k\left[ y_1,\ldots,y_m\right]$.

 Consider an arbitrary Lie polynomial $f(x) \in \rad_F
(Y)$ and a point $$p=\alpha_1 v_1 + \cdots + \alpha_m v_m +c\in
\V.$$ We treat the coefficients $\alpha_i$'s as variables and
write $f(p)$ in the form $$f(p)=g_1(\alpha_1, \ldots,
\alpha_m)u_1+\cdots+ g_s(\alpha_1, \ldots, \alpha_m)u_s,$$ where
$g_1,\ldots,g_s \in k\left[ y_1,\ldots,y_m\right]$ and
$u_1,\ldots,u_s \in F$ are linearly independent and do not depend
on the point $p$. Clearly
\[
f(p)=0 \hbox{ if and only if } g_1(\alpha_1, \ldots, \alpha_m)=
\ldots=g_s(\alpha_1, \ldots, \alpha_m)=0
\]
We therefore can set $S_f=\left\{g_1,\ldots,g_s\right\} \subset
k\left[ y_1,\ldots,y_m\right]$.

 We next construct the inverse correspondence. Consider
a polynomial
\[
g(y_1 ,\ldots,y_m ) = \sum\limits_{\bar i} {\alpha _{\bar i} } y_1
^{i_1 }  \cdots  y_m ^{i_m }  = 0,\quad \alpha _{\bar i} \in
k,\quad \bar i = (i_1 ,\ldots,i_m ) \in \mathbb{N}^m
\]
and set $M_1  = \max \left\{ i_1 \right\}, \ldots, M_m  = \max
\left\{ i_m \right\}$. We associate with the polynomial $g(\bar
y)$ a Lie polynomial $f_g(x)$ such that $$g(\bar y) \in \rad(Y_k)
\leftrightarrow f_g(x)\in\rad(Y_F).$$

Set
$$
\begin{array}{ccc}
f_m (x) & = & s_{m - 1} (x-c,v_1 ,\ldots,v_{n - 1} )  \\
f_{m - 1} (x) & = & s_{m - 1} (x-c,v_1 ,\ldots, v_{n - 2} ,v_n )  \\
   & \vdots &  \\
f_1 (x) & = & s_{n -1} (x-c,v_2 ,\ldots,v_{n - 1} ,v_n )  \\
\end{array},
$$
where  $s_{n - 1}$ is the Lie polynomial defined by Equation
(\ref{eq:III-s_k}) in Section \ref{ss:III-findimalgset}.
Substituting  the point $p$ into the $f_i$'s we obtain
\[
f_m (p) = \alpha _m b_m; \quad f_{m - 1} (p) = \alpha _{m - 1}
b_{m - 1}; \quad \ldots;\quad f_1 (p) = \alpha _1 b_1,
\]
where $b_1,\dots, b_m\in F$ are not zeroes. Choose a Lie
polynomial $a\in F$ so that the degree (see \cite{Bah} for
definition) of $a$ is greater than the degrees of the polynomials
$b_i$'s. In that case all the products of the form $a b_1 \cdots
b_1 b_2 \cdots b_2 \cdots b_n \cdots b_n $ are nonzero (see
\cite{Shirshov} or \cite{BK}). We define the Lie polynomial $f_g
(x)$ as
\begin{gather} \notag
\begin{split}
f_g (x) = \sum\limits_{\bar i} {\alpha _{\bar i} } a \circ
\underbrace {f_1 (x) \circ \cdots \circ f_1 (x)}_{i_1 } \circ
\underbrace {b_1  \circ \cdots \circ b_1 }_{M_1  - i_1 } & \circ
\cdots  \\
\cdots \circ \underbrace {f_m (x) \circ \cdots \circ f_m
(x)}_{i_m} & \circ \underbrace {b_n  \circ \cdots \circ b_n }_{M_m
- i_m }.
\end{split}
\end{gather}
We get, after substituting the point $x=p$ into $f_g (x)$,
\[
f_g(p)=g(\alpha _1 ,\ldots,\alpha _n )a \circ \underbrace {b_1
\circ \cdots \circ b_1 }_{M_1 } \circ \cdots \circ \underbrace
{b_m \circ \cdots \circ b_m }_{M_m }.
\]
Therefore, $f_g (x) = 0$ if and only if $g(\alpha _1 ,\dots
,\alpha _n ) = 0$.

We now define the correspondences $$\rad(Y_F ) \rightarrow
\rad(S_k)$$ and $$\rad(Y_k ) \rightarrow \rad(S_F)$$ by setting
$$
\rad(Y_F) \rightarrow  S_k  = \left\{ S_f \;|\; f \in \rad(Y_F
)\right\}$$ and $$\rad(Y_k)\rightarrow S_F  = \left\{ f_g (x) \in
\Fx \; |\; g \in \rad(Y_k)\right\} \cup s_m (x - c,v_1 ,\ldots,v_m
),
$$
where $V_F(s_m (x - c,v_1 ,\ldots,v_m ))=\V$; see Equation
 (\ref{eq:III-s_k}).

The previous argument holds not only in the case when $Y_F \subset
V+c$ but also in a more general context $Y_F\subseteq F^n$.
Furthermore,
$$\V =V_1+c_1\times \cdots \times V_n+c_n \subset F^n, $$
with $$V_i=\verb"lin"_k\left\{v_1^i, \ldots, v_m^i\right\}.$$
Indeed, similarly to the dimension $1$ case, we can define the
correspondence between algebraic sets from $\V$ and affine subsets
from $k^M$, where $M=m_1+\cdots+m_n$, by setting
\begin{gather} \label{eq:FtoktoFgen}
\begin{split}
    & Y_F=\left\{(\alpha_1^1v_1^1+\cdots+ \alpha_{m_1}^1 v_{m_1}^1,
      \ldots, \alpha_1^nv_1^n+\cdots+ \alpha_{m_1}^n v_{m_1}^n
      +c_n)\right\} \subseteq \V\\  & \leftrightarrow
   Y_k = \{ \underbrace {\alpha_1^1,\ldots,\alpha_{m_1}^1}_{m_1}, \ldots,
    \underbrace {\alpha_1^n,\ldots,\alpha_{m_n}^n}_{m_n}\}
    \subseteq k^M.
\end{split}
\end{gather}

The following theorem demonstrates  relations between this
construction and
correspondence (\ref{eq:FtoktoFgen}).

\begin{thm}\
\begin{itemize}
    \item Let $Y_k \subseteq k^M$ be an algebraic set over $k$.
    Then the corresponding set $Y_F \subset \V$ is algebraic over
    $F$. Moreover, $Y_F=V(S_F)$.
    \item Let $Y_F \subseteq \V$ be an algebraic set over $F$.
    Then the corresponding set $Y_k$ is algebraic over
    $k$. Moreover, $Y_k=V(S_k)$.
    \item The maps $Y_F \rightarrow Y_k$ and $Y_k \rightarrow Y_F$
    define a one-to-one correspondence $Y_F \leftrightarrow Y_k$
    between algebraic {\rm(}over $F${\rm)} sets from
    $\V\subset F$ and algebraic {\rm(}over $k${\rm)} sets from $k^M$.
    Consequently, {\rm $Y_F\rightarrow Y_k \rightarrow
    Y_F=\id_{\AS(F)}$} and {\rm $Y_k\rightarrow Y_F \rightarrow
    Y_k=\id_{\AS(k)}$}.
    \item The correspondence {\rm $\rad(Y_F)\leftrightarrow
    \rad(Y_k)$}
    defined above is a one-to-one correspondence between radical ideals
    from $$k \left[y_1^1,\dots, y_{m_1}^1,\ldots, y_1^n, \dots,y_{m_n}^n \right]$$ and those radical
    ideals from $\F\left[x_1,\dots,x_n\right]$ that contain $$s_{m_1}(x_1-c_1), \ldots, s_{m_n}(x_n-c_n).$$
    Furthermore, to obtain the radical {\rm $\rad(Y_k)$ ($\rad(Y_F)$)} it suffices to
    interpret as algebraic {\rm(}Lie{\rm)} polynomials only those Lie
    {\rm(}algebraic{\rm)}
    polynomials that generate {\rm $\rad(Y_F)$ ($\rad(Y_k)$)}.
\end{itemize}
\end{thm}

\begin{cor}\
\begin{itemize}
    \item Let $S$  be a system of equations over $F$. If \/ $S$
    defines a bounded algebraic set, then $S$ is equivalent to a
    finite system $S_0$.
    \item If the ground field $k$ is finite then any
    subset $M \subset \V$ is algebraic over $F$.
\end{itemize}
\end{cor}

A realisation of a coordinate ring over $k$ can be defined by
analogy with Definition \ref{defn:realisation}. The correspondence
(\ref{eq:FtoktoFgen}) can be reformulated in terms of coordinate
algebras as follows.

\begin{thm}
Let $Y_F \subset \V$ be a bounded algebraic set over $F$ and let
$Y_k \subseteq k^M$ be the corresponding algebraic set over $k$.
Consider the $F$-algebra $$C_F=\left<\x_1,\dots,\x_n\right>_F,$$
where the $\x_i$'s have the form {\rm(\ref{eq:xi})}.

$C_F$ is a realisation of $\Gamma_F(Y_F)$ if and only if the
$k$-subring $$C_k=\left<k,t_1^1,\ldots,t_{m_1}^1, \ldots, t_1^n,
\ldots, t_{m_n}^n \right>$$ is a realisation of the coordinate
ring $\Gamma_k(Y_k)$.
\end{thm}

We have to address some of the differences between the case of
dimension 1 and the case of higher dimensions. If $V = V_1 \times
\cdots \times V_m $ is a fixed variety in $k^N$, we interpret it
in $k^N$ by partitioning the $N$ variables into $m$ groups, $N =
n_1 + \cdots + n_m$. Notice that the affine space $k^N$, $N>1$
represented as sum $m$ sub-spaces is not equivalent to the whole
$k^N$, since the notion of a `decomposed affine space' can not be
expressed in terms of morphisms of algebraic sets.

Notice that the interpretations of bounded algebraic geometry over
the algebra $F$ in dimension 1 and in higher dimensions {\rm(}as
diophantine geometry over $k${\rm)} are very similar.

Finally, we classify  algebraic sets defined by systems of
equations with one variable. By Proposition \ref{prop:boundsA} the
coordinate algebra of any bounded algebraic set has a realisation
in $B(\bar F)$. Moreover, according to Theorem \ref{rem:I-diminB},
$\Gamma (Y)$ is isomorphic to the subalgebra $\left< F,\xi
\right>$ of $B(\bar F)$. Using the machinery of combinatorics in
Lie algebras \cite{BK, Shirshov}, one can show that the only
alternatives are
\begin{itemize}
    \item $\xi \in B(\bar F)$, in which case $Y$ is bounded, or
    \item $\xi \notin B(\bar F)$, in which case $\Gamma (Y)\cong F \ast \langle
    x\rangle$.
\end{itemize}

\begin{thm}
Every algebraic set defined by a system of equations in one
variable over the free Lie algebra $F$ is, up to isomorphism,
either  bounded, or
 empty, or
coincides with $F$.

\end{thm}

Notice that the  Jacoby identity is not essential for the proofs
of most of the results of Section {\rm\ref{sec:AGOFLA}} and they
can be generalised to free anti-commutative algebras
\cite[Appendix]{AGOFLA}.

\centerline{\textsc{Acknowledgements}}

\medskip

The author is grateful to R.~Bryant, E.~Daniyarova and R. St\"or
for useful comments and remarks, as well as to S.~Rees, A.~Duncan,
G.~Megyesi, K.~Goda, M.~Batty for their support and to
H.~Khudaverdyan for lending me his office. My special thanks go to
the London Mathematical Society for the unique opportunity to
write this survey. This text would have never appeared without
V.~N.~Remeselnnikov, to whom I express my sincere gratitude.

\end{document}